\documentclass[12pt]{article}
\usepackage[usenames,dvipsnames]{color}
\usepackage{amsmath,amsthm,amsfonts,amssymb,multicol,amscd,amsbsy}
\usepackage{graphicx}
\usepackage{booktabs}
\usepackage{array}
\usepackage{subfigure}
\usepackage{enumerate}
\usepackage{url}
\usepackage[nodayofweek]{datetime}
\usepackage{mathtools}
\usepackage[toc,page]{appendix}
\usepackage{verbatim} 
\usepackage{amsthm} 
\usepackage{enumitem}
\usepackage{enumerate}
\usepackage{bbm} 
\usepackage[normalem]{ulem} 
\usepackage{bm}
\usepackage{hyperref}
\usepackage{dsfont}

\DeclarePairedDelimiter\ceil{\lceil}{\rceil}
\DeclarePairedDelimiter\floor{\lfloor}{\rfloor}
\usepackage{fullpage}
\usepackage{authblk}

\newcommand{\be}{\begin}
\newcommand{\e}{\end}
\newcommand{\beq}{\begin{equation}}
\newcommand{\eeq}{\end{equation}}

\renewcommand{\l}{\left}
\renewcommand{\r}{\right}

\renewcommand{\d}{\mathrm{d}} 

\newcommand{\ee}{\mathrm{e}}

\newcommand{\set}[1]{\mathbb{#1}}

\newcommand{\curly}[1]{\mathcal{#1}}
\newcommand{\goth}[1]{\mathfrak{#1}}
\newcommand{\setof}[2]{\left\{ #1\; : \;#2 \right\}}

\newcommand{\R}{\set{R}}
\newcommand{\N}{\set{N}}
\newcommand{\C}{\set{C}}
\newcommand{\Z}{\set{Z}}

\newcommand{\E}{\mathbb{E}}

\newcommand{\om}{\omega}
\newcommand{\eps}{\varepsilon}
\newcommand{\Om}{\Omega}

\newcommand{\lam}{\lambda}

\newcommand{\gam}{\gamma}
\newcommand{\Gam}{\Gamma}

\newcommand{\al}{\alpha}
\newcommand{\de}{\delta}

\newcommand{\pr}{\mathbb{P}}		

\newcommand{\ttmatrix}[4]{\left(\be{array}{cc} #1&#2\\	#3&#4 \e{array}	\right)}
\newcommand{\tvector}[2]{\left(\be{array}{c}#1\\#2\e{array}\right)}

\newcommand{\ci}{\mathrm{i}}

\newcommand{\norm}[1]{\left\lVert#1\right\rVert}
\newcommand{\tr}{\mathrm{tr}}	

\newtheorem{thm}{Theorem}[section]

\newtheorem{cor}[thm]{Corollary}
\newtheorem{prop}[thm]{Proposition}

\theoremstyle{definition}
\newtheorem{defn}[thm]{Definition}

\newtheorem{ass}[thm]{Assumption}

\numberwithin{equation}{section}

\theoremstyle{remark}

\def\dotuline{\bgroup
  \ifdim\ULdepth=\maxdimen  
   \settodepth\ULdepth{(j}\advance\ULdepth.4pt\fi
  \markoverwith{\begingroup
  \advance\ULdepth0.08ex
  \lower\ULdepth\hbox{\kern.15em .\kern.1em}%
  \endgroup}\ULon}

\def\dashuline{\bgroup
  \ifdim\ULdepth=\maxdimen  
   \settodepth\ULdepth{(j}\advance\ULdepth.4pt\fi
  \markoverwith{\kern.15em
  \vtop{\kern\ULdepth \hrule width .3em}%
  \kern.15em}\ULon}
\allowdisplaybreaks

\begin{document}
\title{Quantitative lower bounds on the Lyapunov exponent from multivariate matrix inequalities}
\author[1]{Marius Lemm}
\author[2]{David Sutter}
\affil[1]{\small{Department of Mathematics, Harvard University, USA}}
\affil[2]{\small{Institute for Theoretical Physics, ETH Zurich, Switzerland}}
\date{}
\maketitle

\begin{abstract}
The Lyapunov exponent characterizes the asymptotic behavior of long matrix products. Recognizing scenarios where the Lyapunov exponent is strictly positive is a fundamental challenge that is relevant in many applications. In this work we establish a novel tool for this task by deriving a quantitative lower bound on the Lyapunov exponent in terms of a matrix sum which is efficiently computable in ergodic situations. Our approach combines two deep results from matrix analysis --- the $n$-matrix extension of the Golden-Thompson inequality and the Avalanche-Principle. We apply these bounds to the Lyapunov exponents of Schr\"odinger cocycles with certain ergodic potentials of polymer type and arbitrary correlation structure. We also derive related quantitative stability results for the Lyapunov exponent near aligned diagonal matrices and a bound for almost-commuting matrices.
\end{abstract}

\section{Introduction}

Understanding the behavior of long matrix products is a fundamental task that arises in many areas of physics and mathematics, involving the analysis of discrete and continuous dynamical systems and disordered materials~\cite{bougerol_book,carmona_book,viana_book,Wilkinson}. A central question concerning a long matrix product is how fast its norm grows. More precisely, given a sequence $\{L_k\}_{k \geq 1}$ of $d\times d$-dimensional matrices one is interested in the exponential rate
\begin{align} \label{eq_alpha}
\frac{1}{n}\log\l\|\prod_{k=1}^n L_k \r\|,
\end{align}
and its behavior as $n\to\infty$; specifically, one is often interested whether it remains bounded away from zero as $n\to\infty$. Here, $\|\cdot\|$ denotes the operator norm. The quantity \eqref{eq_alpha} is relevant, e.g., in the study of matrices $\{L_k\}_{k \geq 1}$ that are generated by a cocycle over a stationary dynamical system, a very general situation that is commonplace in applications. 
More precisely, let $(X,\mu,T)$ be a dynamical system, i.e., $(X,\mu)$ is a probability space and $T:X\to X$ a measure-preserving dynamical transformation. Given a matrix-valued map $A:X\to \mathrm{GL}_d(\R)$, the dynamics then generate the following cocycle on the trivial bundle over $X$
$$
\begin{aligned}
\curly{C}_A:X\times \R^d &\to X\times \R^d
\\ (x,v)&\mapsto (Tx,A(x) v)  \, .
\end{aligned}
$$
Iterating the cocycle leads to long matrix products taken along orbits of the dynamical transformation $T$, i.e.,
$$
A(T^{n-1} x) A(T^{n-2} x)\ldots A(T x) A(x) \, .
$$
If we additionally assume that $\E \log_+ \| A \|<\infty$ holds with respect to the measure $\mu$,where $(x)_+=\max\{x,0\}$, then the (maximal or top) \textit{Lyapunov exponent} is defined by
\begin{align} \label{eq_topLyap}
\gamma_1:= \lim_{n\to\infty} \frac{1}{n} \E \log \l\|\prod_{k=n-1}^0 A(T^k \cdot) \r\| \, .
\end{align}
The limit in~\eqref{eq_topLyap} exists due to Fekete's subadditivity lemma.

The Lyapunov exponent is especially meaningful in the case where the transformation $T$ is ergodic. In this case, one can $\mu$-almost surely remove the expectation value, i.e., 
\beq\label{eq:kingman}
\begin{aligned} 
\gamma_1
= \lim_{n\to\infty} \frac{1}{n} \log \l\|\prod_{k=1}^n A(T^k \cdot) \r\| \quad  \text{$\mu$-a.s.} \, ,
\end{aligned}
\eeq
by Kingman's subadditive ergodic theorem~\cite{kingman_73} (see also \cite{furstenberg1960}). An important extension is given by Oseledec's theorem \cite{Oseledec}.

A central question about an ergodic cocycle is whether its associated Lyapunov exponent $\gam_1$ is strictly positive or not. (A common situation in applications is that $A:X\to \mathrm{SL}_d(\R)$ in which case $\gam_1\geq 0$ holds automatically.) The strict positivity of the Lyapunov exponent has significant dynamical content and, loosely speaking, indicates chaotic behavior, in the sense of exponential divergence of nearby dynamical trajectories. See \cite{Wilkinson} for a wide-ranging survey on the usefulness of Lyapunov exponents. 

One of the few existing tools for establishing positivity of the Lyapunov exponent is a famous result of Furstenberg~\cite{furstenberg63,furstenberg1971}. Furstenberg's theorem is of immense power and scope, but it also suffers from the disadvantages that it is (a) only applicable to i.i.d.\ matrix products in $\mathrm{SL}_d(\R)$ and (b) entirely non-quantitative. Furstenberg's theorem has been improved and generalized by later works \cite{FurKif,Gor,Kif,Ruelle}, cf.\ the survey \cite{Fu}, in particular by studying the dependence of the Lyapunov exponent on the underlying probability distribution in the i.i.d.\ case. However, the above-mentioned drawbacks (a) and (b) have not been completely removed. Alternative tools for deriving upper and lower bounds on the Lyapunov exponent in special circumstances are provided in \cite{protasov13,sutter19}.

The goal of this paper is to present a new method for deriving lower bounds on the Lyapunov exponent $\gam_1$ based on recently developed techniques from matrix analysis. In fact, we will prove lower bounds on the finite product $\frac{1}{n}\log\l\|\prod_{k=1}^n L_k \r\|$ for every finite $n$ which behave well in the $n\to\infty$ limit in which the Lyapunov exponent emerges, especially in ergodic situations. Our approach combines two powerful and relatively new tools from matrix analysis --- one originating in quantum information theory and the other originating in the spectral theory of Schr\"odinger operators:

\be{enumerate}
\item A recent multivariate trace inequality \textit{generalizing the well-known Golden-Thompson inequality to arbitrarily many matrices} \cite{SBT16}. The original motivation for the result were applications to entropy inequalities in quantum information theory~\cite{Sutter_book}. \label{it_1}
\item
The \emph{Avalanche Principle} which describes the norm of long products of matrices whose expanding directions are somewhat aligned. The Avalanche Principle was originally developed by Goldstein-Schlag~\cite{GS}, for studying the transfer matrices of one-dimensional discrete Schr\"odinger operators. Here we use and refine a recent rendition with effective constants~\cite{HLS} for real-valued matrices that is based on the projective geometry approach of Duarte-Klein \cite{DK1,DK2}.\label{it_2}
\e{enumerate}

It turns out that these two tools combine surprisingly effectively and yield quantitative lower bounds on $\frac{1}{n}\log\l\|\prod_{k=1}^n L_k \r\|$ under perturbatively stable assumptions on the underlying matrix sequence $\{L_k\}_{k \in [n]}$. The bounds require essentially no algebraic structure of the matrix cocycles, though algebraic information may improve the resulting bound.  

The paper is organized as follows. In \textbf{Section \ref{sec_multivariate}}, we review the key tools~\ref{it_1}.\ and~\ref{it_2}.\ mentioned above. In \textbf{Section \ref{sect:main}}, we state our main general results: Theorem \ref{thm:main}, which is the version for finite $n$, and Corollary \ref{cor:asymptotic2}, which is the limiting case in which the Lyapunov exponent arises. In \textbf{Section \ref{sect:schrodinger}}, we present some applications of these new bounds to transfer matrices of one-dimensional discrete ergodic Schr\"odinger operators. These applications yield quantitative lower bounds for the Lyapunov exponent for polymer-type models. Their strength is that they apply to quite general ergodic environments with essentially arbitrary correlation structure (e.g., Ising-correlated random variables and skew-shift dynamics) which we believe is beyond the reach of existing methods.\footnote{For experts on Schr\"odinger operators, we mention that we cannot identify the spectrum very precisely for these models and so our results do not definitively establish positivity of the Lyapunov exponent on the spectrum. Still the bounds we derive compare favorably to the Combes-Thomas estimates \cite{CT,K} even if all the energies where they hold happen to lie outside of the spectrum. See Section \ref{sect:schrodinger} for more details.}\textbf{ Section \ref{sect:pert}} gives two quantitative stability results for the Lyapunov exponent near sufficiently aligned sequences. While these follow from the effective Avalanche Principle alone, they are in a similar spirit as the results in Section \ref{sect:main} in that they apply to general sequences of matrices satisfying preturbatively stable, quantitative and non-algebraic assumptions. Finally, in \textbf{Section \ref{sec_almostComm}}, we describe our related efforts towards proving a lower bound on the Lyapunov exponent that is perturbatively stable for almost-commuting matrices, with the details deferred to the appendix, and we state an open problem which is of interest to applications.

Overall, we would like to summarize our efforts here by saying that the combination of tools~\ref{it_1}.\ and~\ref{it_2}.\ yields coarse, but effective bounds for establishing positivity of the Lyapunov exponent whenever one is in a situation where the effective Avalanche Principle can be verified, which essentially means that the underlying transfer matrices have some degree of hyperbolicity. We certainly hope that the bounds in Theorem \ref{thm:main} and Corollary \ref{cor:asymptotic2} will find further applications in other dynamical contexts. In particular, we wish to point out that all applications considered here only use two scales and utilize an implicit or explicit parameter to ensure the validity of the effective Avalanche Principle deterministically. We believe that these techniques can become even more powerful within an inductive multiscale scheme using also large deviation estimates where the validity of the Avalanche Principle needs to be verified only for typical sequences at the previous scale. 

Another possibility would be to extend the effective Avalanche Principle to complex-valued matrices which would allow, for example, to apply these methods to the complexified Lyapunov exponents appearing in Avila's global theory of one-frequency Schr\"odinger operators \cite{Avila}. In this way, they could potentially be used to distinguish the critical and subcritical phases, thereby distinguishing the absolutely continuous versus singular spectrum. We hope to explore these ideas further in future work. 

\section{Multivariate matrix inequalities} \label{sec_multivariate}
In this section we recall two multivariate matrix inequalities which we eventually combine to bound the Lyapunov exponent from below.
\subsection{Multivariate Golden-Thompson inequality}
The celebrated Golden-Thompson inequality~\cite{golden65,thompson65} states that for any two Hermitian matrices $H_1$ and $H_2$ we have
\begin{align} \label{eq_2GT}
\tr\, \ee^{H_1 +H_2} \leq \tr\, \ee^{H_1} \ee^{H_2} \, .
\end{align}
The inequality \eqref{eq_2GT} was extended to three matrices \cite{Lieb73} and recently to arbitrarily many matrices~\cite{SBT16}, which will be relevant for our results in this manuscript. The result applies to any Schatten $p$-norm $\|\cdot\|_p$ with $1\leq p\leq \infty$, though, we will only use it for $p=\infty$ (with $\|\cdot\|_\infty=\|\cdot\|$ the operator norm).

\be{thm}[$n$-matrix Golden-Thompson~\cite{SBT16}]\label{thm_nGT}
Let $\{H_{k}\}_{k\in[n]}$ be a sequence of Hermitian matrices for $n \in \N$ and $p\geq 1$. Then
\begin{align} \label{eq:nGT}
\log \l\| \exp \left( \sum_{k=1}^n H_k \right) \r\|_p \leq 
\int_{-\infty}^\infty  f(t) \log \l\| \prod_{k=1}^n  \exp\bigl( (1+\ci t) H_k \bigr) \r\|_p \d t \ ,
\end{align}
where $f$ is a probability density on $\R$ defined by
\begin{align} \label{eq_beta0}
f(t) = \frac{\pi}{2} \left( \cosh(\pi t) +1 \right)^{-1} \, .
\end{align}
\e{thm}
For $p=2$, $n=3$, and substituting $H_k \to \frac{1}{2}H_k$ we obtain
\begin{align*}
\tr\,\ee^{H_1 + H_2 + H_3} \leq \int_{-\infty}^\infty  f(t) \tr\, \ee^{H_1} \ee^{\frac{1+\ci t}{2}H_2} \ee^{H_3}\ee^{\frac{1-\ci t}{2}H_2}   \d t \, ,
\end{align*}
which coincides with Lieb's triple matrix inequality~\cite{Lieb73,SBT16}.
For $n=2$ we recover the original Golden-Thompson inequality~\eqref{eq_2GT}. 
We note that there are several related multivariate matrix inequalities which we however not discuss here. The interested reader can find further information in~\cite[Chapter~4]{Sutter_book}.

The main point about the inequality \eqref{eq:nGT} is that it bounds the norm of a matrix product from below by an expression involving the sum of matrices. Sums of matrices are much easier to analyze than products. On a high level, this is because they are commutative quantities. On a technical level, the ergodic theorem describes the asymptotics of normalized sums well and so \eqref{eq:nGT} is a potentially useful tool for analyzing the top Lyapunov exponent. 

\subsection{The Avalanche Principle}
The Avalanche Principle was originally proved by Goldstein-Schlag \cite{GS} for $\mathrm{SL}_2(\R)$-matrices and subsequently generalized to $\mathrm{GL}_d(\C)$-matrices by Schlag \cite{Schlag}. These original versions only applied to sequences of bounded length. More recently, Duarte-Klein \cite{DK1,DK2} gave a new proof based on projective geometry which applies to real-valued matrix sequences of arbitrary length. Here we employ a recent effective rendering of the Duarte-Klein method, where explicit constants were obtained \cite{HLS}.
\begin{defn}
Given a sequence $\{L_k \}_{k \in [n]}$ of matrices in $\mathrm{GL}_d(\R)$, their \emph{expansion rift} is the ratio
\begin{align*}
\rho(L_1, L_2,\ldots,L_{n}):=\frac{\|L_{n}\cdots L_2 L_1\|}{\|L_{n}\|\cdots \|L_2\| \|L_1\|}\in (0,1] \, .
\end{align*}
\end{defn}

Given $L\in \mathrm{GL}_d(\R)$, let $s_1(L)\geq s_2(L)\geq \ldots \geq s_d(L)>0$ denote the sorted singular values of $L$. The largest singular value $s_1(L)$ is the operator norm, i.e.,
$$s_1(L)=\max_{x\in \R^d\setminus \{0\}}\frac{\|Lx\|}{\|x\|}=:\|L\| \, .$$
The smallest singular value of $L$ is the least expansion factor of $L$, regarded as a linear transformation, and it can be characterized by 
$$s_d(L)=\min_{x\in \R^d\setminus \{0\}}\frac{\|Lx\|}{\|x\|}=\|L^{-1}\|^{-1}.$$

\be{defn}
The \emph{gap ratio} (or the \emph{singular gap ratio}) of $L\in\mathrm{GL}_d(\R)$ is the ratio between its first and second singular values, i.e.,
$$\mathrm{gr}(L):=\frac{s_1(L)}{s_2(L)}\, .$$
\e{defn}

The Avalanche Principle is a tool that describes $\|L_n\ldots L_1\|$ in terms of norms of pairs $\|L_{j+1} L_j\|$ in cases where (a) the most expanding directions are significantly more expanding than the other directions (which is quantified by a large gap ratio) and (b) the most expanding directions of subsequent matrices are uniformly somewhat aligned.

Since we utilize an effective Avalanche Principle, it is convenient to introduce notation that allows for various constants to change. We denote $[n]=\{1,2,\ldots,n\}$.

\be{defn}[Effective Avalanche Principle]\label{defn:AP}
Let $\eps_0,c_0,c_l,c_u>0$. We say $$\mathbf{AP}(\eps_0,c_0,c_l,c_u)$$ if the following holds: For every $n\in \N$, $0<\eps \leq \eps_0$ and $0<\kappa\leq c_0\eps^2$, and every finite sequence $\{L_k\}_{k \in [n]}$ of matrices in $\mathrm{GL}_d(\R)$, satisfying the two assumptions
\begin{enumerate}
\item[(G)] \textit{Gap:} $\quad\mathrm{gr}(L_i)\geq \kappa^{-1}$ for all $i \in [n]$ \label{ass_G}
\item[(A)]  \textit{Alignment:} $\quad\rho(L_{i}, L_{i+1})\geq \eps$ for all $i \in [n-1]$\, , \label{ass_A}
\end{enumerate}
it holds that
\begin{align*}
\ee^{-c_l n\kappa/{\eps^2}}
\leq \frac{\rho(L_1,L_2,\ldots,L_{n})}{\rho(L_1,L_2)\cdots \rho(L_{n-1}, L_{n})}
\leq \ee^{c_u n\kappa/{\eps^2}}.
\end{align*}
\e{defn}

The following version of an effective Avalanche Principle was proved in Section 5 of \cite{HLS}.

\be{thm}[Effective AP version 1 \cite{HLS}]\label{thm:AP1}
We have $\mathbf{AP}\l(\frac{1}{10},\frac{1}{10},5,11\r)$.
\e{thm}

In the derivation of Theorem \ref{thm:AP1} in \cite{HLS} various parameter choices were made which were appropriate for the purposes considered there but which can be modified with an eye towards other applications, e.g., the ones we consider in Section \ref{sect:schrodinger} here. We state a relaxed version of the effective Avalanche Principle with weaker assumptions on $\kappa$ and $\eps_0$ (i.e., with larger constants $\eps_0$ and $c_0$).

\be{thm}[Effective AP version 2]\label{thm:AP2}
We have $\mathbf{AP}\l(\frac{1}{5},\frac{1}{6},11,11\r)$ for all matrix products of length $n\geq 36$.
\e{thm}

The error constants $c_l=c_u=11$ are comparable to the ones in Theorem \ref{thm:AP1} above. The restriction to $n\geq 36$ is unimportant for applications. The proof of Theorem \ref{thm:AP2} follows the general line of argumentation in \cite{HLS}, which is an effective version of \cite{DK1,DK2}, but is subject to various modifications and refinements that we describe in Appendix~\ref{ap_AP}. 

The main reason why we include Theorem \ref{thm:AP2} is as follows: We will see below in \eqref{eq:APworstcase} that the validity of the assumptions (G) and (A) of $\mathbf{AP}\l(\frac{1}{10},\frac{1}{10},5,11\r)$ for all matrices along a sequence implies its uniform hyperbolicity (meaning that the Lyapunov exponent is positive irrespective of the underlying dynamics). This is not the case for $\mathbf{AP}\l(\frac{1}{5},\frac{1}{6},11,11\r)$ with its weaker assumptions. This is relevant because for Schr\"odinger operators, uniform hyperbolicity holds exactly on the complement of the spectrum by Johnson's theorem \cite{Johnson,Zhang} and so it is useful to have version 2 of the effective Avalanche Principle which can apply on the spectrum as well.




\be{rmk}
The statement of the Avalanche Principle in \cite{DK1,DK2,HLS} also contains information on the most expanding directions (the eigenvectors associated to $s_1$) for the long product  $\prod_{k=n}^1 L_k$. We omit that part of the statement since we will not use it in the following.
\e{rmk}

\section{The main result}\label{sect:main}
In this section, we present our main general result, a quantitative lower bound on the Lyapunov exponent for a product of invertible matrices satisfying the effective Avalanche Principle.

\subsection{The bound for finite products}
Our first main result is the following lower bound on the finite product $\frac{1}{n}\log\l\|\prod_{k=1}^n L_k \r\| $. It applies to sequences of invertible matrices satisfying an effective Avalanche Principle. We write $\lambda_{\max}(\cdot)$ for the maximal eigenvalue of a matrix and $|A|=\sqrt{A^\dagger A}$ for its absolute value.  


\be{thm}[Bounds for finite products]\label{thm:main}
For $\eps_0,c_0,c_l,c_u>0$, assume that $\mathbf{AP}(\eps_0,c_0,c_l,c_u)$ holds. Moreover, assume that $\{L_k\}_{k \in [n]}$, $n \in \N$, are matrices in $\mathrm{GL}_d(\R)$ that satisfy the assumptions (G) and (A) of $\mathbf{AP}(\eps_0,c_0,c_l,c_u)$ with constants $\kappa,\eps>0$.
\be{enumerate}[label=(\roman*)]
\item  \label{it_main_i} Suppose that each $L_k$ is a normal matrix. Then
\beq\label{eq:thmmain1}
\frac{1}{n}\log\l\|\prod_{k=1}^n L_k \r\| \geq \lam_{\max}\l(\frac{1}{n}\sum_{k=1}^n \log |L_k|\r)-\frac{(c_l+c_u)\kappa}{\eps^2} \, .
\eeq
\item  \label{it_main_ii}
Suppose that there exists $\al\in [0,1)$ such that
\beq\label{eq:strange}
\|L_{k+1} L_k\|\geq (1-\al)\| |L_{k+1}^\dagger| |L_k|\| \qquad \forall j\in [n]\, .
\eeq
Then
\beq\label{eq:thmmain2}
\begin{aligned}
 \frac{1}{n}\log\l\|\prod_{k=1}^n L_k \r\| \geq&
\lam_{\max} \l(\frac{1}{n}\sum_{k=1}^n \frac{\log |L_k|+\log |L_k^\dagger|}{2}\r)
 -\l(c_l+\frac{c_u}{(1-\al)^2}\r)\frac{\kappa}{\eps^2}+\log(1-\al).
 \end{aligned}
\eeq
\e{enumerate}
\e{thm}

\be{rmk}\label{rmk:main}\
\be{enumerate}[label=(\alph*)]
\item
Statement~\ref{it_main_ii} is a relaxation of the case of normal matrices considered in statement~\ref{it_main_i} which shows that the estimate is robust. Indeed, for normal matrices, Lemma \ref{lm:norms} below implies that the assumption \eqref{eq:strange} holds with $\al=0$ and using $|L_k^\dagger|=|L_k|$ one sees that the estimate \eqref{eq:thmmain2} with $\al=0$ reduces to \eqref{eq:thmmain1}.

\item In the commutative case where $[L_k,L_{k'}]=[L_k,L_{k'}^\dagger]=0$ for all $k,k'\in [n]$, we find
\beq\label{eq_Lyap_commute}
\begin{aligned} 
 \frac{1}{n}\log\l\|\prod_{k=1}^n L_k \r\|
&= \frac{1}{n} \log \norm{\prod_{k=1}^n |L_k|}
= \frac{1}{n} \log \norm{\exp\Big( \sum_{k=1}^n \log |L_k| \Big)}\\ 
 &=\lambda_{\max} \left( \frac{1}{n} \sum_{k=1}^n \log |L_k| \right) \, 
\end{aligned}
\eeq
so \eqref{eq:thmmain1} can be understood as saying that $ \frac{1}{n}\log\l\|\prod_{k=1}^n L_k \r\|$ is bounded from below by its value in the commutative case up to errors coming from the Avalanche Principle. This makes sense insofar as the idea that commutators can be ignored if one is interested in particular inequalities is precisely the core message of the Golden-Thompson inequality and its generalizations.

\item The important class of Schr\"odinger cocycles where each $L_k$ is a product of matrices of the form $[a,-1;1,0]$ with $a\in \R$ falls under the scope of statement~\ref{it_main_ii}. Indeed, \eqref{eq:strange} holds with $\al=0$ in the Schr\"odinger case. This follows from the observation that $U_* [a,-1;1,0] U_*^{-1}=A^\dagger$ with the $a$-independent unitary $U_*=[1,0;0,-1]$; see Lemma \ref{lm:strange}.
\e{enumerate}
\e{rmk}

Theorem \ref{thm:main} is proved in Section \ref{sect:mainproof} via the effective Avalanche Principle and the $n$-matrix Golden-Thompson inequality.

As promised in Remark \ref{rmk:main}(a), we verify that condition \eqref{eq:strange} holds with $\al=0$ for normal matrices.

 \be{lm}[Norm identities]\label{lm:norms}\
\be{enumerate}[label=(\roman*)]
\item
Let $A$ and $B$ be arbitrary matrices. Then $\|AB\|=\||A| |B^\dagger|\|$.
\item If $A$ and $B$ be are normal normal matrices, then also $\|AB\|=\||A^\dagger| |B|\|$.
\e{enumerate}
\e{lm}

\be{proof}
Let $Y$ be an arbitrary matrix. The identities $\|Y\|^2=\|Y^\dagger Y\|=\|YY^\dagger\|=\|Y^\dagger\|^2$ imply 
$$
\|AB\|^2=\|B^\dagger |A|^2 B\|=\||A| |B^\dagger|^2 |A|\|=\||A| |B^\dagger|\|^2 \,,
$$
for arbitrary $A$ and $B$. If $A$ and $B$ are normal, then $|A|=|A^\dagger|$ and $|B^\dagger|=|B|$. 
\e{proof}

%

\subsection{The lower bound on Lyapunov exponents}
We apply Theorem \ref{thm:main} to derive a quantitative lower bound on the Lyapunov exponent of ergodic matrix cocycles. This follows immediately from the finitary Theorem \ref{thm:main}.

\be{cor}\label{cor:asymptotic2}
Let $\{L_k\}_{k \in \N}$ be a cocycle of matrices in $\mathrm{GL}_d(\R)$ over a an ergodic dynamical system $(X,\mu,T)$. Assume that $\mathbb E\log \|L_1\|<\infty$.
\be{enumerate}[label=(\roman*)]
\item Suppose that the set of $x\in X$ such that the assumptions of Theorem \ref{thm:main}(i) are satisfied for the entire sequence $\{L_k\}_{k \in \N}$ has positive $\mu$-measure. Then
\beq\label{eq:cor1}
\gam_1
\geq \lam_{\max}\l(  \mathbb E \log |L_1|\r)-\frac{(c_l+c_u)\kappa}{\eps^2} \, .
\eeq

\item Suppose that the set of $x\in X$ such that the assumptions of Theorem \ref{thm:main}(ii) are satisfied for the entire sequence $\{L_k\}_{k \in \N}$, has positive $\mu$-measure. Then
\beq\label{eq:cor2}
\gam_1
\geq \lam_{\max}\l( \frac{\mathbb E \log |L_1|+  \mathbb E \log |L_1^\dagger|}{2}\r)-\l(c_l+\frac{c_u}{(1-\al)^2}\r)\frac{\kappa}{\eps^2}+\log(1-\al) \, .
\eeq
\e{enumerate}
\e{cor}

\be{proof}
This follows directly from Theorem \ref{thm:main}, Kingman's subadditive ergodic theorem in the form \eqref{eq:kingman}, the standard ergodic theorem, and continuity of the function $\lam_{\max}(\cdot)$.
\e{proof}

These bounds are powerful in the presence of a large, but fixed, parameter that allows to verify the validity of the effective Avalanche Principle. Under such circumstances, they yield quantitative lower bounds on Lyapunov exponents which are also perturbatively stable. To avoid confusion, we note that the Avalanche Principle $\mathbf{AP}(\eps_0,c_0,c_l,c_u)$ with constants $\kappa,\eps>0$ for the $\{L_k\}_{k\in\N}$ by itself already implies an a priori lower bound on the Lyapunov exponent
\beq\label{eq:APworstcase}
\begin{aligned}
\gam_1&\geq \liminf_{n\to\infty}\frac{1}{n}\sum_{k=1}^{n-1}\log\frac{\|L_{k+1}L_k\|}{\|L_k\|} -\frac{c_l\kappa}{\eps^2}
\geq \liminf_{n\to\infty}\frac{1}{n}\sum_{k=1}^{n}\log(\eps\|L_k\|)-\frac{c_l\kappa}{\eps^2}\\
&\geq \frac{1}{2}\log\l(\frac{\eps^2}{\kappa}\r)-\frac{c_l\kappa}{\eps^2},
\end{aligned}
\eeq
Since the Avalanche Principle is a deterministic statement, this is a worst-case type bound. Still, it may yield a positive lower bound on $\gam_1$ in the case when $\frac{\kappa}{\eps^2}$ and $c_l$ are sufficiently small, and in those cases it implies uniform hyperbolicity of all cocycles ranging over the set of $\{L_k\}$. This is the case, e.g., for Theorem \ref{thm:AP1} where it is assumed that $\frac{\kappa}{\eps^2}\leq c_0=\frac{1}{10}$ but not for the Theorem \ref{thm:AP2} which has weaker assumptions.

 Generally speaking, Corollary~\ref{cor:asymptotic2} is useful in situations where the average case, as represented by $\lam_{\max}\l( \frac{\mathbb E \log |L_1|+  \mathbb E \log |L_1^\dagger|}{2}\r)$, behaves better than the worst-case bound \eqref{eq:APworstcase}, and we give a prototypical example of this in the following subsection.

\subsection{A prototypical example}
We illustrate the usefulness of Corollary~\ref{cor:asymptotic2} with a brief example. Additional applications to Schr\"odinger cocycles are postponed to the next section. 

Consider the shift $T$ on the space of two-sided sequences $X=\{0,1\}^\Z$. Let $\mu_0$ be an ergodic measure on $X$ with
\beq\label{eq:gothpdefn0}
\goth{p}:=\mu_0\l(\setof{\om\in X}{\om_0=0}\r).
\eeq
Furthermore, let $ \sqrt{1000}\leq a\leq b$ and let $R(\theta)=(\cos \theta, -\sin \theta; \sin \theta, \cos \theta)$ be the rotation matrix with $\theta \in [0,\pi]$. We define a cocycle $A:X\to \mathrm{SL}_2(\R)$ by 
\begin{align*}
A(0):=\begin{pmatrix} a & 0 \\ 0 & a^{-1} \end{pmatrix} ,\qquad
A(1):= R(\pi/4) \begin{pmatrix} b & 0 \\ 0 & b^{-1} \end{pmatrix} R(\pi/4)^\dagger 
\end{align*}

We denote $A(i)=A_i$ for $i\in\{0,1\}$. Elementary considerations show that $\mathrm{gr}(A_0)  = a^2$, $\mathrm{gr}(A_1)  = b^2$, $\rho(A_0,A_0)=\rho(A_1,A_1)=1$ and
$
\rho(A_0,A_1)\geq \frac{1}{\sqrt{2}}.
$
(The last inequality is an asymptotic equality for large $a$- and $b$-values.) Hence, for any $x\in X=\{0,1\}^\Z$, the sequence of normal matrices $L_k:=A(T^k x)$, $k\in\N$, induced by the cocycle satisfies the assumptions of $\mathbf{AP}\l(\frac{1}{10},\frac{1}{10},5,11\r)$ with $\kappa=a^{-2}\leq \frac{1}{1000}$ and $\eps=\frac{1}{10}$. 

First, the worst-case bound \eqref{eq:APworstcase} then reads
\beq\label{eq:ex1}
\gam_1\geq\frac{1}{2}\log\l(\frac{a^2}{100}\r)-\frac{500}{a^2}=\log a-\log 10-\frac{500}{a^2}.
\eeq
Note \eqref{eq:ex1} is independent of the parameters $b\geq a$ and the probability $\goth{p}$ from \eqref{eq:gothpdefn0}. Especially the lack of $b$-dependence means that this bound, even if it gives a positive number, does not capture the size of the Lyapunov exponent. To be more precise, we note the following upper bound\footnote{We remark that a tighter upper bound in the form of a convex optimization problem that can be evaluated efficiently for any fixed $a$ and $b$ is presented in~\cite{sutter19}.} on $\gam_1$, which holds by submultiplicativity of the norm,
\begin{align}\label{eq:ex3}
\gamma_1 \leq \E \log \norm{L_1} = \goth{p} \log a + (1-\goth{p}) \log b \, .
\end{align}
Note the discrepancy between the bounds \eqref{eq:ex1} and \eqref{eq:ex3}.

 Now, the lower bound \eqref{eq:ex1} can be substantially improved with the help of Corollary~\ref{cor:asymptotic2} (i) which gives 
\beq\label{eq:ex2}
\begin{aligned}
\gamma_1 
&\geq \lam_{\max}\l(  \mathbb E \log |L_1|\r)-\frac{1600}{a^2}
=\lam_{\max}\l(  \goth{p} \log |A_0|+(1-\goth{p})\log |A_1|\r)\\
&= \sqrt{(1-\goth{p})^2(\log b)^2+\goth{p}^2(\log a)^2}-\frac{1600}{a^2},
\end{aligned}
\eeq
Indeed, the bound \eqref{eq:ex2} manages to capture the size of the Lyapunov exponent much more accurately than the worst-case bound \eqref{eq:ex1} whenever $\goth{p}\neq 1$ and $b>a$;  compare with the upper bound \eqref{eq:ex3}. We emphasize that the bound \eqref{eq:ex2} only depends on the underlying dynamics through the parameter $\goth{p}$ and consequently applies for arbitrarily correlated ergodic dynamical systems.

To see that the bound \eqref{eq:ex2} is quite accurate, we note that in the limit where $a,b\to\infty$ with $\frac{a}{b}\to 0$ (and assuming $\goth{p}< 1$), the upper bound \eqref{eq:ex3} and the lower bound \eqref{eq:ex2} match asymptotically and together imply the precise asymptotic $\gam_1\sim (1-\goth{p})\log b$. 

\subsection{Proof of the main result, Theorem \ref{thm:main}}\label{sect:mainproof}
To prove Theorem \ref{thm:main} we combine the effective Avalanche Principle with the $n$-matrix Golden-Thompson inequality. 
Remarkably, these two rather deep results on long matrix products work hand-in-hand: The $n$-matrix Golden-Thompson inequality \eqref{eq:nGT} has the structural advantage of bounding a long matrix product in terms of a long matrix sum. However, it suffers from the technical disadvantage that the matrix product is dressed with lots of $t$-dependent unitaries which also get averaged in $t$. The key observation is that the Avalanche Principle perfectly cures this technical ailment because it allows us to replace the long matrix product by the product over pairwise terms $\| L_{j+1}L_j\|$, and precisely in this setting the unitaries do not play a significant role (at least under the assumptions in (i) or (ii) of Theorem \ref{thm:main}) thanks to unitary invariance of the operator norm.

\be{proof}[Proof of Theorem \ref{thm:main}]
By the assumption that $\mathbf{AP}(\eps_0,c_0,c_l,c_u)$ holds and that $\{L_k\}_{k \in [n]}$, $n \in \N$, satisfy its assumptions, we can use the lower bound in the Avalanche Principle to obtain
\beq\label{eq:thmstep1}
 \al_n=\frac{1}{n}\log\l\|\prod_{k=1}^n L_k\r\|
 \geq  \frac{1}{n}\sum_{j=1}^{n-1}\log\l(\rho(L_j,L_{j+1})\r) +\frac{1}{n}\sum_{j=1}^{n}\log\|L_j\| -\frac{c_l\kappa}{\eps^2} \, .
\eeq
By Lemma \ref{lm:norms}, statement (i) follows from statement (ii), so we assume that we are in case (ii) of Theorem \ref{thm:main}.  
 
  Fix an arbitrary $t\in \R$. The fact that $\|A\|=\||A|\|=\||A^\dagger|\|$ for arbitrary matrices combined with unitary invariance of the operator norm implies that
\beq\label{eq:id1}
 \|L_j\|= \||L_j|\|=\||L_j|^{1+it}\|=\||L_j^\dagger|^{1+it}\|.
\eeq
Using this identity, Assumption \eqref{eq:strange}, and unitary invariance of the norm, we obtain
\beq\label{eq:id2}
\begin{aligned}
\rho(L_j,L_{j+1})
&=\frac{\|L_{j+1}L_j \|}{ \|L_{j+1}\|\| L_j\|}\\
&\ge (1-\al)\frac{\| |L_{j+1}^\dagger| |L_j| \|}{ \||L_{j+1}^\dagger|\|\| |L_j|\|}\\
&= (1-\al)\frac{\| |L_{j+1}^\dagger|^{1+\ci t}  |L_j|^{1+\ci t}  \|}{ \||L_{j+1}^\dagger|^{1+\ci t}\|\| |L_j|^{1+\ci t}\|}\\
&= (1-\al) \rho\l(|L_j|^{1+\ci t},|L_{j+1}^\dagger|^{1+\ci t}\r)\, .
\end{aligned}
\eeq
Similarly, replacing Assumption \eqref{eq:strange} by an application of Lemma \ref{lm:norms} (i), we have
\beq\label{eq:id3}
\rho(L_j,L_{j+1})=\frac{\| |L_{j+1}| |L_j^\dagger| \|}{ \||L_{j+1}|\|\| |L_j^\dagger|\|}=\frac{\| |L_{j+1}|^{1+\ci t}  |L_j^\dagger|^{1+\ci t}  \|}{ \||L_{j+1}|^{1+\ci t} \|\| |L_j^\dagger|^{1+\ci t} \|}=\rho\l(|L_j^\dagger|^{1+\ci t},|L_{j+1}|^{1+\ci t}\r).
\eeq

Using \eqref{eq:id2} and \eqref{eq:id3} we can estimate the right-hand side of \eqref{eq:thmstep1} in terms of
the alternating sequence
$$
\tilde L_j:=
\be{cases}
L_j,\textnormal{ if $j$ is odd},\\
L_j^\dagger,\textnormal{ if $j$ is even},
\e{cases}
$$
 as follows
$$
\begin{aligned}
&\frac{1}{n}\sum_{j=1}^{n-1}\log\l(\rho(L_j,L_{j+1})\r) +\frac{1}{n}\sum_{j=1}^{n}\log\|L_j\| -\frac{c_l\kappa}{\eps^2}\\
&\hspace{20mm}\geq \frac{1}{n}\sum_{j=1}^{n-1}\log\l(\rho(|\tilde L_j|^{1+\ci t},|\tilde L_j|^{1+\ci t})\r) +\frac{1}{n}\sum_{j=1}^{n}\log\||\tilde L_j|^{1+\ci t}\| -\frac{c_l\kappa}{\eps^2}+\log(1-\al)
\end{aligned}
$$
Next, observe that
$$
\mathrm{gr}(L_j)=\mathrm{gr}\l(|\tilde L_j|\r)=\mathrm{gr}\l(|\tilde L_j|^{1+it}\r) \, .
$$
Together with \eqref{eq:id2} and \eqref{eq:id3}, this implies that the matrices $\{|\tilde L_k|^{1+\ci t}\}_{k\in [n]}$ inherit the validity of assumptions (G) and (A) in $\mathbf{AP}(\eps_0,c_0,c_l,c_u)$ from the $\{L_k\}_{k\in [n]}$ but with $\eps$ replaced by $(1-\al)\eps$.

Hence, we can apply the upper bound in the Avalanche Principle and obtain
$$
\begin{aligned}
 \al_n
 \geq&  \frac{1}{n}\sum_{j=1}^{n-1}\log\l(\rho(|\tilde L_j|^{1+\ci t},|\tilde L_j|^{1+\ci t})\r) +\frac{1}{n}\sum_{j=1}^{n}\log\||\tilde L_j|^{1+\ci t}\| -\frac{c_l\kappa}{\eps^2}+\log(1-\al)\\
 \geq&   \frac{1}{n}\log\l\|\prod_{k=1}^n |\tilde L_k|^{1+\ci t}\r\| -\l(c_l+\frac{c_u}{(1-\al)^2}\r)\frac{\kappa}{\eps^2}+\log(1-\al),
  \end{aligned}
$$
for every fixed $t\in \R$. 

We can average both sides of this inequality over $t$ with respect to the probability measure $ f(t)\d t$ defined in~\eqref{eq_beta0}. Afterwards, we are in a position to apply the $n$-matrix Golden-Thompson inequality \eqref{eq:nGT} from Theorem~\ref{thm_nGT}. This gives
\beq\label{eq:thmstep2}
\begin{aligned}
  \al_n 
 &\geq   \frac{1}{n} \int_{\mathbb R}  f(t)\log\l\|\prod_{k=1}^n |\tilde L_k|^{1+\ci t}\r\| \d t-\l(c_l+\frac{c_u}{(1-\al)^2}\r)\frac{\kappa}{\eps^2}+\log(1-\al)\\
 &\geq \frac{1}{n}\log\l\|\exp\l(\sum_{k=1}^n \log |\tilde L_k|\r)\r\| -\l(c_l+\frac{c_u}{(1-\al)^2}\r)\frac{\kappa}{\eps^2}+\log(1-\al)\\
 &= \lam_{\max}\l(\frac{1}{n}\sum_{k=1}^n \log |\tilde L_k|\r) -\l(c_l+\frac{c_u}{(1-\al)^2}\r)\frac{\kappa}{\eps^2}+\log(1-\al).
 \end{aligned}
\eeq
In the last step, we used that the spectral theorem for $X=\sum_{k=1}^n \log |\tilde L_k|$ implies that $\|\exp(X)\|=\exp(\lam_{\max}(X))$.

Finally, we repeat the argument with the sequence $\{\tilde L_k\} $ replaced by the other alternating sequence
$$
L_j':=\tilde L_j^\dagger=
\be{cases}
L_j^\dagger,\textnormal{ if $j$ is odd},\\
L_j,\textnormal{ if $j$ is even},
\e{cases}
$$
invoking again the identities \eqref{eq:id1}--\eqref{eq:id3}. This gives
\beq\label{eq:thmstep3}
\al_n\geq  \lam_{\max}\l(\frac{1}{n}\sum_{k=1}^n \log |L_k'|\r) -\l(c_l+\frac{c_u}{(1-\al)^2}\r)\frac{\kappa}{\eps^2}+\log(1-\al).
\eeq
Taking the average of the estimates \eqref{eq:thmstep2} and \eqref{eq:thmstep3} and using subadditivity of $\lam_{\max}$ yields
$$
\al_n\geq \lam_{\max}\l(\frac{1}{n}\sum_{k=1}^n \frac{\log |L_k|+\log |L_k^\dagger|}{2}\r)-\l(c_l+\frac{c_u}{(1-\al)^2}\r)\frac{\kappa}{\eps^2}+\log(1-\al),
$$
and thus the assertion of Theorem \ref{thm:main} in case (ii). 
\end{proof}

\section{Applications to ergodic Schr\"odinger cocycles}
\label{sect:schrodinger}

In this section, we apply the main results to Schr\"odinger cocycles of polymer type. More precisely, we consider potentials  $v$ which take fixed values for $2p$ steps before a new value is ergodically sampled. Here $p\in \N$ is a fixed integer and we are interested in establishing a positive Lyapunov exponent. The random dimer case ($p=1$) was first studied by Dunlap-Wu-Philips \cite{DWP} who observed delocalization at special energies; a thorough mathematical investigation of this phenomenon was then conducted in \cite{JSB,JSBS}.  For us, the reason for studying polymers is purely technical: For a sufficiently large polymer length parameter $p$, we are able to verify the validity of the effective Avalanche Principle for a single polymer block and this yields efficiently computable positive lower bounds on the Lyapunov exponent. These bounds only require little information on the underlying ergodic dynamics, essentially just the individual probability of occurrence for each polymer type (and there are just two types in the simplest case). Consequently, the bounds apply to ergodic dynamical systems with arbitrary correlation structures, which we believe is not achievable by other methods in this generality. We mention that some of the most difficult problems in the analysis of dynamically defined Schr\"odinger cocycles involve the case of ``intermediate correlations'' (i.e., not too rigid and not too random), such as, e.g., skew-shift or doubling map dynamics. All of these dynamical systems are covered by Theorem \ref{thm:schro} below.

The main question regarding Schr\"odinger cocycles is whether the Lyapunov exponent is positive inside the spectrum of the associated discrete Schr\"odinger operator on $\ell(\Z)$. This positivity of the Lyapunov exponent is closely related to the electronic phenomenon of Anderson localization. Off the spectrum, the Lyapunov exponent is a priori known to be strictly positive by the well-known Combes-Thomas estimates \cite{CT,K}. In fact, Johnson's result \cite{Johnson} characterizes the spectrum as the complement of the set where the cocycle is uniformly hyperbolic. 

Unfortunately, in the polymer models we consider, we cannot locate the spectrum precisely and we therefore cannot conclude positivity of the Lyapunov exponent inside the spectrum of any particular model. Still, as a consolation prize, we can say that there is spectrum very close to the points where we obtain a lower bound on the Lyapunov exponent. This means that even if our bounds only apply at energies that happen to lie outside of the spectrum, they are considerably stronger than the Combes-Thomas estimates which explicitly deteriorate as one approaches the spectrum.



\subsection{Ergodic Schr\"odinger cocycles}
To define the cocycles, we let $(X,\mu,T)$ be an ergodic dynamical system and let $f:X\to \R$ be measurable and bounded. For any energy parameter $E\in \R$, we define the cocycle by the map $A_E:X\to \mathrm{GL}_d(\R)$,
\beq\label{eq:schrodefn}
A_E(x):=\ttmatrix{E-f(x)}{-1}{1}{0},\qquad x\in X.
\eeq
These arise as the transfer matrices of the discrete Schr\"odinger operator $H(x):\ell^2(\Z)\to\ell^2(\Z)$ defined by
\beq\label{eq:Hdefn}
(H(x)\psi)_n=\psi_{n+1}+\psi_{n-1}+v_n(x)\psi_n\qquad \textnormal{with }  v_n(x):=f(T^n x)\qquad \forall n\in \Z\, .
\eeq
Such $H(x)$ are called Schr\"odinger operators with dynamically defined potentials $v_n$, for general background, see the recent survey \cite{Damanik}. The precise connection to \eqref{eq:schrodefn} is that if $\{u_n\}_{n\in \Z}$ solves $Hu=Eu$, then it can be expressed only in terms of its initial data via the transfer matrices, i.e.,
$$
\tvector{u_{n}}{u_{n-1}}=\prod_{j=n-1}^0 A(T^j x) \tvector{u_{1}}{u_0} \qquad \forall n\in \Z \, .
$$
The operator family $\{H(x)\}_{x\in X}$ is a covariant ergodic family in the sense that
\beq\label{eq:Hergodic}
H(Tx)=S H(x) S^{-1},
\eeq
where $(S\psi)_n=\psi_{n+1}$ is the left shift, a unitary on $\ell^2(\Z)$. We recall a fundamental result of Pastur for such ergodic operator families.

\be{thm}[Non-randomness of the spectrum, Pastur \cite{Pastur}]\label{thm:pastur}
There exists a set $\Sigma\subset \R$ such that
$$
\mathrm{spec}\, H(x)=\Sigma,\qquad \mu-a.e.\ x\in X.
$$
\e{thm}

%
%

\subsection{Formal definition of the polymer models}
We fix a parameter $p\in\N$, called the polymer length. The potential sequence $\{v_n\}_{n\in \Z}$ of the polymer Schr\"odinger operator only takes two values $\{0,-v\}\in \R$ for some constant $v>0$. From these, we construct the polymer blocks of length $2p$
\beq\label{eq:hatvdefn}
\hat v_+=(0,\ldots,0,-v,\ldots,-v)\in \R^{2p},
\qquad \hat v_-=(-v,-v,\ldots,-v)\in \R^{2p},
\eeq
where the number of $0$s in $\hat v_+$ is equal to $p$.

In a nutshell, in the polymer model the two blocks $\hat v_+$ and $\hat v_-$ are sampled according to a probability measure $\mu_0$ on the space $\Om=\{+1,-1\}^\Z$ which we assume is ergodic with respect to the left shift $\curly{S}$ on $\Om$.

We now give the formal construction of the polymer operators that are ergodic in the sense of \eqref{eq:Hergodic}. This is slightly cumbersome for the technical reason described in Remark \ref{rmk:technical} (i) below and can be skipped on a first reading. We follow the formalism in \cite{JSB,JSBS}, where the i.i.d.\ case is considered. We begin by decomposing $\Om$ into the sets
$$
\Om=\Om_+\cup \Om_-,\qquad \textnormal{where } \Om_\pm:=\setof{\om \in \Om}{\om_0=\pm}.
$$
Then the total probability space on which we define our model is
\beq\label{eq:Xdefn}
 X:=X_+\cup X_-, \qquad \textnormal{where } X_\pm:=\Om_\pm\times \{0,1,\ldots,p-1\},
\eeq
endowed with the probability measure
\beq\label{eq:mudefn}
\mu(A_\pm\times \{\ell\}):=\frac{\mu_0(A_\pm )}{p},\qquad \forall A_\pm \subset \Om_\pm,\quad \forall \ell\in \{0,1,\ldots,p-1\}.
\eeq
We define the invertible transformation $T:X\to X$ by
\beq\label{eq:Tdefn}
T(\om,\ell)=
\be{cases}
(\om,\ell+1), &\textnormal{if } \ell<p-1\\
(\curly{S}\om,0), &\textnormal{if } \ell=p-1\, ,
\e{cases}
\eeq
where $\curly{S}$ is the left shift on $\Om$.

\be{prop}
$(X,\mu,T)$ defined by \eqref{eq:Xdefn}--\eqref{eq:Tdefn} is an ergodic dynamical system.
\e{prop}

\be{proof}
This follows from the assumption that $\curly{S}:\Om\to\Om$ is an ergodic transformation with respect to the measure $\mu_0$. We omit the details.
\e{proof}

We are now ready to give a formal definition of polymer potential.

\be{defn}[Formal definition of polymer potential]
\label{defn:poly}
Let $p\in \N$, and let $v_+\in \R$. Define the constant vectors $\hat v_+,\hat v_-\in \R^{2p}$ by \eqref{eq:hatvdefn}. Define the measurable and bounded sampling function $f:X\to\R $ by
$$
f(\om,l):=\hat v_{\om_0}(\ell) \, .
$$
This defines a family of Schr\"odinger operators $\{H(\om,\ell)\}_{(\om,\ell)\in X}$ via \eqref{eq:Hdefn} with $x=(\om,\ell)$. 
\e{defn}

We call this the family of Schr\"odinger operators with polymer length $2p$ and polymer potentials $\hat v_+,\hat v_-\in\R^{2p}$ sampled according to $\mu_0$. It is easy to verify that it is a covariant ergodic family in the sense of \eqref{eq:Hergodic}. This completes the formal construction.

\be{rmk}\label{rmk:technical}
\
\be{enumerate}[label=(\roman*)] 
\item
As mentioned in the beginning, at the heart of the construction is the probability measure $\mu_0$ which describes the sampling of the polymer sequence. The remainder of the construction involving $X,\mu,T$ only serves to generate a Schr\"odinger operator which is covariant under the standard left-shift as in \eqref{eq:Hergodic}. The construction achieves this by appropriately randomizing the location of the initial polymer. We mention that a simpler alternative is to always fix the first site of the first polymer to be $0\in\Z$ and then to sample polymers according to $\mu_0$. This simpler alternative yields a family of Schr\"odinger operators which is still ergodic, however not with respect to the shift $S$ by a single site but with respect to the shift $S^{2p}$ by $2p$ sites, i.e., from polymer to polymer. We prefer the slightly more technical setup here since it leads to a family of operators which is ergodic in the standard way expressed by \eqref{eq:Hergodic}.

\item
The setup and methods straightforwardly generalize to other potentials \eqref{eq:hatvdefn} and to the case where the hopping coefficients along the polymers are also sampled ergodically, i.e., when the ergodic Schr\"odinger operator is replaced by an ergodic Jacobi operator of polymer type; see, e.g., \cite{JSB}.
\e{enumerate}
\e{rmk}

\subsection{Spectral information for the polymer models}
We identify points that are close to the spectrum of the polymer models. We write $\goth{p}$ for the probability that $\mu_0$ has the polymer $\hat v_-$ as its first piece, i.e.,
\beq\label{eq:gothpdefn}
\goth{p}:=\mu_0(\setof{\om\in \Om}{\om_0=-})\in [0,1].
\eeq

The following proposition shows that the polymer models have at least some spectrum in $[-2,2]$ and the number of spectral points grows with $p$. Here $\mathrm{dist}$ denotes the distance on $\R$.

\be{prop}[Spectral information]\label{prop:spec}
Suppose that $\goth{p}<1$ and let $k\in[p]$. Then  
$$
\mathrm{dist}\l(2\cos\l( \frac{\pi k}{p+1} \r),\,\mathrm{spec}\, H(x)\r)
\leq 18\frac{k}{p^{3/2}},
$$
holds for $\mu_0$-a.e. $x\in X$.
\e{prop}
\be{proof}
Let $k\in[p]$ and $\al=\frac{\pi k}{p+1}$. Define the $\ell^2(\Z)$-normalized sequence $\phi=(\phi_n)_{n\in\Z}$ by
$$
\phi_n=
\be{cases}
\sqrt{\frac{2}{p}}\sin(\al n),\qquad &\textnormal{if } 1\leq n\leq p,\\
0,\qquad &\textnormal{otherwise}.
\e{cases}
$$
By the assumption $\mu_0(\Om_-)>0$ and Theorem \ref{thm:pastur}, we can assume that $x=(\om,\ell)\in X$ satisfies $\om_0=-$ and $\ell=0$. Then, the facts that $\phi_{n+1}+\phi_{n-1}=2\cos(\al)\phi_n$ for $2\leq n\leq p-1$ and $\sin x\leq x$ imply
\beq\label{eq:weyl}
\begin{aligned}
\|\l(H(x)-2\cos\al\r)\phi\|
&= \sqrt{\phi_1^2+\phi_p^2+(\phi_2-2\phi_1\cos\al)^2+(\phi_{p-1}-2\phi_p\cos\al)^2}\\
&\leq \sqrt{5\phi_1^2+5\phi_p^2+2\phi_2^2+2\phi_{p-1}^2}
=\sqrt{10\phi_1^2+4\phi_2^2}\\
&=\sqrt{\frac{2}{p}}\sqrt{10\sin^2\al+4\sin^2(2\al)}
\leq \sqrt{\frac{32}{p}}\al=\sqrt{\frac{32}{p}}\frac{\pi k}{p+1}\\
&\leq 18\frac{k}{p^{3/2}}.
\end{aligned}
\eeq
Since $\|\phi\|=1$, we can define a constant Weyl sequence that is always equal to $\phi\in \ell^2(\Z)$. The bound \eqref{eq:weyl} then proves Proposition \ref{prop:spec}.
\e{proof}

\subsection{Lower bound on the Lyapunov exponent in polymer models}

The main result of Section \ref{sect:schrodinger} reads as follows.

\be{thm}[Main result about polymer models]\label{thm:schro}
Let $\goth{p}\geq 1/2$, $p\geq 2$, and let $E=2\cos\theta$ with $\theta\in[0,2\pi]$ so that there exist $0<\de_1,\de_2<\frac{\pi}{2}$ satisfying
\beq\label{eq:Edelta12ass}
 \min\{\theta,2\pi-\theta\}\geq \de_1 
\qquad \text{and} \qquad
 \min_{k\in[p]}\l|(p+1)\theta-k\pi \r|\geq\de_2\, .
\eeq
Assume that $v\geq b_0+E$ with
\beq\label{eq:bass}
b_0:=1+\max\l\{4,\,
\l(\frac{20}{9}\r)^8\frac{1}{\de_1\de_2},
\,\l(\frac{20}{9}\r)^8\frac{1}{(\de_1\de_2)^2}
,\,\l(10^{6}\de_2^{-1}\r)^{2/(2p-1)}
,\, \l(\frac{160}{9\de_1}\r)^{10/(p-1)}
\r\}.
\eeq
Then
\beq
\gam_1\geq \frac{\goth{p}}{2}p\log b>0\, .
\eeq
\e{thm}

The significance of Theorem \ref{thm:schro} is to establish positivity of the Lyapunov exponent for all energies $E$ satisfying \eqref{eq:Edelta12ass} independently of the correlation structure between the $\hat v_+,\hat v_-$ and whether $E$ lies in the spectrum or not. 

\be{rmk}\
\be{enumerate}[label=(\roman*)]
\item The first assumption \eqref{eq:Edelta12ass} says that the energy $E=2\cos\theta$ does lie in the spectrum $[-2,2]$ of the free Laplacian but it is bounded away from the special points that already appeared in Proposition \ref{prop:spec}. The fact that these points are asymptotically spaced at distance $1/p$ for large $p$ inside of $[-2,2]$ implies that even if $E$ is not in the spectrum, there is definitely spectrum near $E$ (in a way that is made precise by Proposition \ref{prop:spec}).

\item Heuristically, the second assumption $v\geq b_0+E$ with $b_0$ relatively large as specified in \eqref{eq:bass} says that the potential value $-v$ of the polymers is rather unfavorable to the energy $E$. 

\item The set of energies $E\in [-2,2]$ that is covered by Theorem \ref{thm:schro} does not contain the spectral points found in Proposition \ref{prop:spec}, but they can lie close to them (they can lie arbitrarily close to these spectral points for sufficiently large $p$, in fact). Hence, even if all the energies $E$ covered by Theorem \ref{thm:schro} happen to lie outside of the spectrum, the bound $\gam_1\geq \frac{\goth{p}}{2}p\log b>0$ proved here is generally stronger than the Combes-Thomas estimates \cite{CT}, due to the two facts that (a) the Combes-Thomas estimates only give a small lower bound on the Lyapunov exponent near the spectrum and (b) there are many (on the order of $p$ for large $p$) spectral points spread out over $[-2,2]$ by Proposition \ref{prop:spec}.
\e{enumerate}
\e{rmk}

\subsection{Proof of Theorem \ref{thm:schro}}

\subsubsection{Norm bounds on matrix powers}

We prove the following bounds on the singular gap and alignment of the $p$th powers of matrices 
\beq\label{eq:AB}
A=\ttmatrix{a}{-1}{1}{0}, \qquad B=\ttmatrix{b}{-1}{1}{0}
\eeq
with $a,b\in\R\setminus\{\pm2\}$. We note that their eigenvalues are given by
\beq\label{eq:lammudefn}
\lam_\pm=\frac{a\pm \sqrt{a^2-4}}{2}\in\C,\qquad \mu_\pm=\frac{b\pm \sqrt{b^2-4}}{2}\in\C,
\eeq
and we denote $\lam=\lam_+$ and $\mu=\mu_+$.

\be{prop}[Norm bounds on matrix powers]\label{prop:APverify}
Let $p,q\in\N$. Let $a=2\cos\theta$ with $\theta\in[0,2\pi]$ and let $0<\de_1,\de_2<\frac{\pi}{2}$ be so that \eqref{eq:Edelta12ass} holds. Moreover, let $b\geq b_0$ with $b_0$ as in \eqref{eq:bass}. Then
\begin{align}
\frac{9}{10}\mu^{q}\leq&\|B^q\|\leq \frac{20}{9}\mu^{q},\\
\frac{9\de_2}{40}\mu^q\leq& \|A^pB^q\|\leq \frac{160}{9\de_1}\mu^{q},\\
\frac{\de_2}{4}\l(\frac{9}{10}\r)^2 \mu^{3p}\leq& \|B^{2p}A^pB^p\|,\\
\frac{1}{2}\l(\frac{9\de_2}{20}\r)^{2}\mu^{2p}\leq& \|A^pB^pA^pB^p\|.
\end{align}
\e{prop}

The following $F$-function governs the growth of these norms.

\be{defn}[$F$-function]\label{defn:F}
Let $q\in \N$, and define the function $F_q:\C\setminus\{0,\pm 1\}\to\C$ by
\beq\label{eq:Fqdefn}
F_q(z):=\frac{z^q-z^{-q}}{z-z^{-1}}
\eeq
where the branch cut for the complex logarithm is placed along the negative real axis.
\e{defn}

We remark that the branch cut will play no important role. 

\be{lm}
\label{lm:pnormbounds}
Let $p,q\in \N$, $a,b\in\R\setminus\{\pm2\}$. Define the matrices $A$ and $B$ by \eqref{eq:AB} and let $\lam=\lam_+$ and $\mu=\mu_+$ be as in \eqref{eq:lammudefn}. Then
\begin{align}
\label{eq:Anorm}
\max_{p'\in \{p-1,p,p+1\}} |F_{p'}(\lam)|\leq &
\|A^p\|\leq 2\max_{q\in \{p-1,p,p+1\}} |F_{p'}(\lam)|,\\
\label{eq:Bnorm}
\max_{q'\in \{q-1,q,q+1\}} |F_{q'}(\mu)|\leq &
\|B^q\|\leq 2\max_{q'\in \{q-1,q,q+1\}} |F_{q'}(\mu)|,\\
 \label{eq:ABnorm}
  |F_{p+1}(\lam)F_{q+1}(\mu)\!-\!F_{p}(\lam)F_{q}(\mu)|\leq &\|A^pB^q\|\leq 4\!\!\max_{p'\in \{p-1,p,p+1\}}\!\! |F_{p'}(\lam)|\max_{q'\in \{q-1,q,q+1\}}|F_{q'}(\mu)|,
\end{align}
and the lower bounds
\begin{align}
   \label{eq:BABnorm}
&\|B^{2p} A^p B^p\|\\
\nonumber
&\geq|F_{2p+1}(\mu)(F_{p+1}(\lam)F_{p+1}(\mu)\!-\!F_{p}(\lam)F_{p}(\mu))
\!-\!F_{2p}(\mu)(F_{p}(\lam)F_{q+1}(\mu)\!-\!F_{p-1}(\lam)F_{q}(\mu))|,\\
    \label{eq:ABABnorm}
  &\|A^{p} B^{q} A^{p} B^{q}\|\\
  \nonumber
  &\geq |(F_{p+1}(\lam)F_{q+1}(\mu)-F_{p}(\lam)F_{q}(\mu))
(F_{p+1}(\lam)F_{q+1}(\mu)-F_{p}(\lam)F_{q}(\mu))\\
\nonumber
&\quad+(-F_{p+1}(\lam)F_{q}(\mu)-F_{p}(\lam)F_{q-1}(\mu))
(F_{p}(\lam)F_{q+1}(\mu)-F_{p-1}(\lam)F_{q}(\mu))|
\end{align}
\e{lm}

\be{proof}[Proof of Lemma \ref{lm:pnormbounds}]
Since $|a|\neq 2$, we can diagonalize $A$ as follows.
\beq\label{eq:Adiag}
\begin{aligned}
A=S D S^{-1},\qquad \textnormal{where }
&D=\ttmatrix{\lam_-}{0}{0}{\lam_+}
\qquad
S=\ttmatrix{\lam_-}{\lam_+}{1}{1},\\
&S^{-1}=\frac{1}{\lam_+-\lam_-} \ttmatrix{-1}{\lam_+}{1}{-\lam_-},
\end{aligned}
\eeq
and $\lam_\pm$ as in \eqref{eq:lammudefn}. Let $p\in\N$. The diagonalization allows to efficiently compute $A^p$. Recalling that $\lam_+\lam_-=1$ and $\lam=\lam_+$, as well as Definition \ref{defn:F}, we have
\beq\label{eq:A^p}
A^p\!=\!S D^p S^{-1}
\!=\!\frac{1}{\lam_+-\lam_-}\ttmatrix{\lam_+^{p+1}-\lam_-^{p+1}}{\lam_-^{p}-\lam_+^{p}}
{\lam_+^{p}-\lam_-^{p}}{\lam_-^{p-1}-\lam_+^{p-1}}
\!=\!\ttmatrix{F_{p+1}(\lam)}{-F_{p}(\lam)}{F_{p}(\lam)}{-F_{p-1}(\lam)}.
\eeq
The estimate \eqref{eq:Anorm} now follows from the matrix norm equivalence
$$
\max_{i,j=1,2}|(A^p)_{i,j}|\leq \|A^p\|\leq 2\max_{i,j=1,2}|(A^p)_{i,j}|,
$$
after recalling Definition \eqref{eq:Fqdefn} of $F_p$.  The matrix $B$ can also be diagonalized, $B=\tilde S\tilde D \tilde S^{-1}$,  with $\tilde S$ and $\tilde D$ defined as in \eqref{eq:Adiag} but with $\lam_\pm$ replaced by $\mu_\pm$ from \eqref{eq:lammudefn}. Hence, the bounds in \eqref{eq:Bnorm} hold by the same arguments.

We come to \eqref{eq:ABnorm} next. The identity \eqref{eq:A^p} and its analog for $B^q$ give
$$
A^p B^q
=\ttmatrix{F_{p+1}(\lam)F_{q+1}(\mu)-F_{p}(\lam)F_{q}(\mu)}{-F_{p+1}(\lam)F_{q}(\mu)+F_{p}(\lam)F_{q-1}(\mu)}{F_{p}(\lam)F_{q+1}(\mu)-F_{p-1}(\lam)F_{q}(\mu)}{-F_{p}(\lam)F_{q}(\mu)+F_{p-1}(\lam)F_{q-1}(\mu)} \, .
$$
To prove the upper bound in \eqref{eq:ABnorm}, we estimate again by the maximum-entry norm followed by the triangle inequality to obtain
$$
\ \|A^p B^q\|\leq 2\max_{i,j=1,2}|(A^p B^q)_{i,j}|\leq 
4\max_{p'\in \{p-1,p,p+1\}} F_{p'}(\lam)\max_{q'\in \{q-1,q,q+1\}}F_{q'}(\mu),
$$
as desired. Similarly, 
$$
\begin{aligned}
 \|A^p B^q\|\geq
\max_{i,j=1,2}|(A^p B^q)_{i,j}|
\geq |(A^p B^q)_{1,1}|
\geq|F_{p+1}(\lam)F_{q+1}(\mu)-F_{p}(\lam)F_{q}(\mu)|
\end{aligned}
$$
and this proves \eqref{eq:ABnorm}. For \eqref{eq:BABnorm}, we lower bound the norm by the top left entry to find
$$
\begin{aligned}
 &\|B^{2p} A^p B^p\|\geq
\max_{i,j=1,2}|(B^{2p} A^p B^p)_{i,j}|
\geq |(B^{2p} A^p B^p)_{1,1}|\\
&\geq|F_{2p+1}(\mu)(F_{p+1}(\lam)F_{p+1}(\mu)-F_{p}(\lam)F_{p}(\mu))
-F_{2p}(\mu)(F_{p}(\lam)F_{q+1}(\mu)-F_{p-1}(\lam)F_{p}(\mu))|.
\end{aligned}
$$
Finally, we again use the lower bound on the norm in terms of the top left entry to find
$$
\begin{aligned}
\|A^{p} B^{q} A^{p} B^{q}\|
&\geq |(A^{p} B^{q} A^{p} B^{q})_{1,1}|\\
&\geq |(F_{p+1}(\lam)F_{q+1}(\mu)-F_{p}(\lam)F_{q}(\mu))
(F_{p+1}(\lam)F_{q+1}(\mu)-F_{p}(\lam)F_{q}(\mu))\\
&\quad+(-F_{p+1}(\lam)F_{q}(\mu)-F_{p}(\lam)F_{q-1}(\mu))
(F_{p}(\lam)F_{q+1}(\mu)-F_{p-1}(\lam)F_{q}(\mu))|
,
\end{aligned}
$$
which proves \eqref{eq:ABABnorm} and hence Lemma  \ref{lm:pnormbounds}.
\e{proof}

We need to control the $F$-function in two regimes: when $z=e^{i\theta}$ with $\theta$ away from a few special points and when $z=r>1$.

\be{lm}[Bounds on the $F$-function]\label{lm:Fbounds}
\
\be{enumerate}[label=(\roman*)] 
\item Let $p\in\N$ and $0<\de_1,\de_2<\frac{\pi}{2}$. For every $z=e^{i\theta}$ with $\theta\in [0,2\pi]$, if $\de_1\leq \min\{\theta,2\pi-\theta\}$, then $|F_p(z)|\leq \frac{2}{\de_1}$, and if $\de_2\leq \min_{k\in[p]}\l|(p+1)\theta-k\pi \r|$ then $|F_{p+1}(z)|\geq \frac{\de_2}{2}$.

\item Let $q\in\N.$ For every $z=x\in\R$, if $x>x_0>1$, then $ x^{q-1}(1-x_0^{-2q})
\leq
F_q(x)
\leq x^{q-1}(1-x_0^{-2})^{-1}
$.
\e{enumerate}
\e{lm}

\be{proof}[Proof of Lemma \ref{lm:Fbounds}]
For part (i), we first assume $\de_1\leq \min\{\theta,2\pi-\theta\}$. By Definition \eqref{eq:Fqdefn} of $F$ and the fact that $\sin x\geq \frac{x}{2}$ for $x\in[0,\pi/2]$ to obtain
$$
|F_q(z)|=|F_q(e^{i\theta)}|=\l|\frac{\sin(q\theta)}{\sin(\theta)}\r|
=\frac{|\sin(q\theta)|}{\sin|\theta|}
\leq \frac{1}{\sin \de_1}\leq \frac{2}{\de_1}.
$$
Second, we assume that $\de_2\leq \min_{k\in[p]}\l|(p+1)\theta-k\pi \r|$ and use similar estimates to obtain
$$
F_{p+1}(z)=\l|\frac{\sin((p+1)\theta)}{\sin(\theta)}\r|
\geq |\sin((p+1)\theta)|
\geq \sin \de_2
\geq \frac{\de_2}{2}.
$$
This proves part (i). For part (ii), we note that $x>x_0>1$ implies $x^{-1}<x_0^{-2}x$. Hence, 
$$
F_q(z)=F_q(x)=\l|\frac{x^q-x^{-q}}{x-x^{-1}}\r|=\frac{x^q-x^{-q}}{x-x^{-1}}
\be{cases}
\leq \frac{x^q}{x-x^{-1}}\leq\frac{x^{q-1}}{1-x_0^{-2}},\\
\geq \frac{x^q-x^{-q}}{x}\geq x^{q-1}(1-x_0^{-2q})
\e{cases}
$$
and Lemma \ref{lm:Fbounds} is proved.
\e{proof}

We are now ready to give the 

\be{proof}[Proof of Proposition \ref{prop:APverify}]
By Lemma \ref{lm:pnormbounds}, it suffices to control the $F$-function in the various cases. We start with $\|B^q\|$. 
 By Definition \eqref{eq:lammudefn} and the fact that $b\geq b_0>2$ by assumption,
$$
\begin{aligned}
\mu=\mu_+= \frac{b+\sqrt{b^2-4}}{2}
\geq b-1=:x_0>1.
\end{aligned}
$$
We note that by \eqref{eq:bass}, the choice $x_0=b-1\geq 4$ satisfies 
\beq\label{eq:x0est}
 1-x_0^{-2p}\geq 1-x_0^{-2}\geq \frac{9}{10}.
\eeq
We apply Lemma \ref{lm:Fbounds} and \eqref{eq:x0est} to obtain
$$
\|B^q\|\geq \max_{q'\in\{q-1,q,q+1\}} F_{q'}(\mu)\geq 
F_{q+1}(\mu)
\geq 
\mu^{p} (1-x_0^{-2p})
\geq\frac{9}{10}\mu^{p}.
$$
and also
\beq\label{eq:Bpub}
\|B^q\|\leq 2\max_{q'\in\{q-1,q,q+1\}} F_{q'}(\mu)
\leq 
2\max_{q'\in\{q-1,q,q+1\}}
\mu_+^{q'-1} (1-x_0^{-2})^{-1}
\leq 2\mu_+^{p} (1-x_0^{-2})^{-1}
\leq \frac{20}{9}\mu^{p}.
\eeq

We come to $\|A^pB^q\|$ next. For this, we combine \eqref{eq:ABnorm} and Lemma \ref{lm:Fbounds}. For the latter, note that $a=2\cos\theta$ is equivalent to $\lam=\frac{a+i\sqrt{4-a^2}}{2}=e^{i\theta}$ and so the assumptions of Proposition \ref{prop:APverify} imply the validity of the conditions for Lemma \ref{lm:Fbounds}(i). Using \eqref{eq:x0est}, we obtain
$$
\|A^pB^q\|
\leq 4\max_{p'\in \{p-1,p,p+1\}} |F_{p'}(\lam)|\max_{q'\in \{q-1,q,q+1\}}F_{q'}(\mu)
\leq \frac{80}{9}\mu^{q}
\max_{p'\in \{p-1,p,p+1\}} |F_{p'}(\lam)|
\leq \frac{160}{9\de_1}\mu^{q}.
$$
and the lower bound
$$
\begin{aligned}
\|A^pB^q\|
&\geq |F_{p+1}(\lam)F_{q+1}(\mu)-F_{p}(\lam)F_{q}(\mu)|\\
&\geq |F_{p+1}(\lam)F_{q+1}(\mu)|-|F_{p}(\lam)F_{q}(\mu)|\\
&\geq \frac{9\de_2}{20} \mu^q-\frac{20}{9\de_1}\mu^{q-1}\\
&= \mu^q \l(\frac{9\de_2}{20} -\frac{20}{9\de_1}\mu^{-1}\r)\\
&\geq \frac{9\de_2}{40}\mu^q,
\end{aligned}
$$
where we used $\mu\geq b-1$ and \eqref{eq:bass} in the last step. 

Next, we combine \eqref{eq:BABnorm}, Lemma \ref{lm:Fbounds} and \eqref{eq:x0est} to obtain the lower bound
$$
\begin{aligned}
\|B^{2p} A^p B^p\|
\geq&
|F_{2p+1}(\mu)(F_{p+1}(\lam)F_{p+1}(\mu)-F_{p}(\lam)F_{p}(\mu))\\
&-F_{2p}(\mu)(F_{p}(\lam)F_{p+1}(\mu)-F_{p-1}(\lam)F_{p}(\mu))|\\
\geq& |F_{2p+1}(\mu)F_{p+1}(\lam)F_{p+1}(\mu)|-|F_{2p+1}(\mu)F_{p}(\lam)F_{p}(\mu))|\\
&-|F_{2p}(\mu)F_{p}(\lam)F_{p+1}(\mu)|-|F_{2p}(\mu)F_{p-1}(\lam)F_{p}(\mu))|\\
\geq& \frac{\de_2}{2}\l(\frac{9}{10}\r)^2\mu^{3p}
-\frac{6}{\de_1}\l(\frac{10}{9}\r)^2 \mu^{3p-1}\\
=& \mu^{3p}\l(\frac{\de_2}{2}\l(\frac{9}{10}\r)^2-\frac{6}{\de_1}\l(\frac{10}{9}\r)^2\mu^{-1}\r)\\
\geq&\frac{\de_2}{4}\l(\frac{9}{10}\r)^2 \mu^{3p}.
\end{aligned}
$$
where we used $\mu\geq b-1\geq b_0-1$ and \eqref{eq:bass} in the last step. 

Finally, we combine \eqref{eq:ABABnorm} at $q=p$ with Lemma \ref{lm:Fbounds} and \eqref{eq:x0est} to find
$$
\begin{aligned}
  \|A^pB^pA^p B^p\|
  \geq
  &|(F_{p+1}(\lam)F_{p+1}(\mu)-F_{p}(\lam)F_{p}(\mu))^2\\
   &+(-F_{p+1}(\lam)F_{p}(\mu)-F_{p}(\lam)F_{p-1}(\mu))  
(F_{p}(\lam)F_{p+1}(\mu)-F_{p-1}(\lam)F_{p}(\mu))|\\
  \geq
  &|F_{p+1}(\lam)F_{p+1}(\mu)-F_{p}(\lam)F_{p}(\mu)|^2\\
   &-|(F_{p+1}(\lam)F_{p}(\mu)+F_{p}(\lam)F_{p-1}(\mu))  
(F_{p}(\lam)F_{p+1}(\mu)-F_{p-1}(\lam)F_{p}(\mu))|\\
\geq& \l(\frac{9\de_2}{20}\r)^{2}\mu^{2p}-7\l(\frac{20}{9\de_1}\r)^2\mu^{2p-1}\\
=& \mu^{2p}\l(\l(\frac{9\de_2}{20}\r)^{2}-7\l(\frac{20}{9\de_1}\r)^2\mu^{-1}\r)\\
\geq& \frac{1}{2}\l(\frac{9\de_2}{20}\r)^{2}\mu^{2p},
\end{aligned}
$$
where we used $\mu\geq b-1\geq b_0-1$ and \eqref{eq:bass} in the last step. This proves Proposition \ref{prop:APverify}.
\e{proof}

\subsubsection{Verification of the conditions of the Avalanche Principle}
We now use Proposition \ref{prop:APverify} to verify the conditions of the effective Avalanche Principle. We define the matrix sequence
\beq\label{eq:Lkrewrite}
L_k:=\prod_{j=2pk-1}^{2p(k-1)} A_E(T^j x) =\prod_{j=2pk-1}^{2p(k-1)} \ttmatrix{E-f(T^jx)}{-1}{1}{0}.
\eeq
Since ergodicity implies that the Lyapunov exponent arises as the limit for almost every $x\in X$, we may restrict without loss of generality to the event that $x=(\om,0)$ in the following, i.e., that the polymer of type $\om_0\in \{\pm\}$ begins at the site $0\in \Z$.

\be{cor}\label{cor:APverify}
Under the assumptions of Theorem \ref{thm:schro}, the matrix sequence $\{L_k\}_{k\in\N}$ satisfies the conditions of $\mathbf{AP}\l(\frac{1}{5},\frac{1}{6},11,11\r)$ with $\kappa=\l(\frac{5}{\de_2}\r)^2\mu^{-2p}$ and $\eps=10^{-4}\min\{1,(\delta_1\delta_2)^2\}$.
\e{cor}

\be{proof}
We recall that by Definition \ref{defn:poly}, $f(\om,\ell)=\hat v_{\om_0}(\ell)$ is the polymer potential with $\hat v_\pm$ given in \eqref{eq:hatvdefn}. Since we can assume $x=(\om_0,0)$, this means that $f(T^j x)=\hat v_{(T^j\om)_0}(T^j 0)$ in \eqref{eq:Lkrewrite}. Hence, by setting
$$
a:=E,\qquad b:=E+v.
$$
we find that the sequence $\{L_k\}_{k\in\Z}$ takes the two (matrix) values,
$$
L_k=
\be{cases}
A^p B^p,\qquad &\textnormal{if } (T^{2p(k-1)}\om)_0=-,\\ 
B^{2p}, \qquad &\textnormal{if } (T^{2p(k-1)}\om)_0=+
\e{cases}
$$
and Assumptions of Theorem \ref{thm:schro}, specifically \eqref{eq:Edelta12ass} and \eqref{eq:bass}, ensure that the conditions of Proposition \ref{prop:APverify} are verified which then implies that Conditions (G) and (A) are verified with $\kappa$ and $\epsilon$ given as follows. First, Proposition  \ref{prop:APverify} gives 
$$
\l(\min\l\{\|A^pB^p\|,\|B^{2p}\|\r\}\r)^{-2}
\!\leq\! \l(\min\l\{\frac{9\de_2}{40}\mu^{p},\, \frac{9}{10}\mu^{2p}\r\}\r)^{-2}
\!\leq\! \l(\frac{40}{9\de_2}\r)^2 \mu^{-2p}
\!\leq\! \l(\frac{5}{\de_2}\r)^2\mu^{-2p}=\kappa,
$$
where the last step uses that $\mu\geq b_0-1\geq 4\geq\frac{\de_2}{4}$ by \eqref{eq:bass}. Second, Proposition \ref{prop:APverify} gives
$$
\begin{aligned}
&\min\l\{\frac{\|A^pB^pA^pB^p\|}{\|A^pB^p\|^2},\frac{\|B^{4p}\|}{\|B^{2p}\|^2},\frac{\|A^pB^{3p}\|}{\|A^pB^p\|\|B^{2p}\|},\frac{\|B^{2p}A^pB^p\|}{\|B^{2p}\|\|A^pB^p\|}\r\}
\geq
\frac{1}{10^4}\min\{1,(\delta_1\delta_2)^2\}
=\eps\, . 
\end{aligned}
$$
This proves Corollary \ref{cor:APverify}.
\e{proof}

\subsubsection{Verification of condition \eqref{eq:strange} with $\al=0$}
In this short subsection, we note an algebraic property of all Schr\"odinger transfer matrices which implies that condition \eqref{eq:strange} holds with $\al=0$, so these matrices behave as good as normal matrices for the purposes of Theorem \ref{thm:main}.

\be{lm}\label{lm:strange}
For $\{L_k\}_{k\in\N}$ as in \eqref{eq:Lkrewrite}, it holds that
$$
\|L_{k+1}L_k\|=\||L_{k+1}^\dagger| |L_k|\|.
$$
\e{lm}

\be{proof}
The key observation is that all Schr\"odinger transfer matrices are unitarily equivalent to their adjoint via the fixed unitary
$$
U_*:=\ttmatrix{1}{0}{0}{-1}=U_*^{-1} \, .
$$
Indeed, direct computation shows that
$$
U_*AU_*=A^\dagger=\ttmatrix{a}{1}{-1}{0}
$$
and since $U_*$ is independent of $a$, we also have $U_*BU_*=B^\dagger$. Hence, invoking unitary invariance of the operator norm and Lemma \ref{lm:norms} (i), 
$$
\|L_{k+1}L_k\|
=\|U_*L_{k+1}U_*^2 L_kU_*\|
=\|L_{k+1}^\dagger L_k^\dagger\|
=\||L_{k+1}^\dagger| |L_k|\|
$$
as claimed. This proves Lemma \ref{lm:strange}.
\e{proof}

\subsubsection{Conclusion}
We are now ready to give the 

\be{proof}[Proof of Theorem \ref{thm:schro}]
As noted above, we may assume that $x=(\om,0)$. By Corollary \ref{cor:APverify}, the matrix sequence $\{L_k\}_{k\in\N}$ from \eqref{eq:Lkrewrite} satisfies the condition of $\mathbf{AP}\l(\frac{1}{5},\frac{1}{6},11,11\r)$ with $\kappa=\l(\frac{5}{\de_2}\r)^2\mu^{-2p}$ and $\eps=10^{-4}\min\{1,(\delta_1\delta_2)^2\}$. By Lemma \ref{lm:strange}, the $\{L_k\}$ matrices satisfy condition \eqref{eq:strange} with $\al=0$. Hence, by applying Corollary \ref{cor:asymptotic2}, recalling \eqref{eq:gothpdefn} and Weyl's theorem on norm perturbations of Hermitian matrices,
\beq\label{eq:conclusion}
\begin{aligned}
\gam_1
\geq& 
\lam_{\max}\l( \frac{\mathbb E \log |L_1|+  \mathbb E \log |L_1^\dagger|}{2}\r)-22\frac{\kappa}{\eps^2}\\
\geq& 
\lam_{\max}\l(\goth{p} \frac{ \log |B^{2p}|+  \log|(B^{2p})^\dagger|}{2}
+(1-\goth{p}) \frac{ \log |A^pB^p|+ \log |(A^pB^p)^\dagger|}{2}\r)\\
&-10^{11}\mu^{-2p}\de_2^{-2}\max\{1,(\de_1 \de_2)^{-2}\}\\
\geq&
\lam_{\max}\l(\goth{p} \frac{ \log |B^{2p}|+  \log|(B^{2p})^\dagger|}{2}\r)\\
&-\frac{1-\goth{p}}{2}\l\|\log |A^pB^p|+ \log |(A^pB^p)^\dagger|\r\|
-10^{11}\mu^{-2p}\de_2^{-2}\max\{1,(\de_1 \de_2)^{-2}\} \, .
\end{aligned}
\eeq
We control the first error term by noting that $\|\log X\|=\log\|X\|$ for  Hermitian matrices $X\in \mathrm{SL}(2,\R)$. Using this fact, $\||X^\dagger|\|=\|X\|$, and Proposition \ref{prop:APverify}, we obtain
$$
\begin{aligned}
\l\|\log |A^pB^p|+ \log |(A^pB^p)^\dagger|\r\|
&\leq \log \||A^pB^p|\|+ \log \||(A^pB^p)^\dagger|\|
= 2\log\|A^pB^p\|\\
&\leq 2\log\l(\frac{160}{9\de_1}\mu^{p}\r)
=2p\log \mu+2\log\l(\frac{160}{9\de_1}\r).
\end{aligned}
$$
Next we compute the main term in \eqref{eq:conclusion}. For this, we note that, after diagonalizing $B$, \eqref{eq:A^p} for $B^{2p}$ reads
$$
B^{2p}=\ttmatrix{F_{p+1}(\mu)}{-F_{p}(\mu)}{F_{p}(\mu)}{-F_{p-1}(\mu)}.
$$
We note that $F_q(\mu)\geq 0$ for $q\in\{2p-1,2p,2p+1\}$ by Lemma \ref{lm:Fbounds} (ii) and \eqref{eq:x0est}. We can employ the operator monotonicity of the logarithm and the square root to estimate
$$
\begin{aligned}
\log|B^{2p}|
&=\log\sqrt{(B^{2p})^\dagger B^{2p}}\\
&=\sqrt{\ttmatrix{F_{2p+1}(\mu)^2+F_{2p}(\mu)^2}{-F_{2p+1}(\mu)F_{2p}(\mu)-F_{2p}(\mu)F_{2p-1}(\mu)}{-F_{2p+1}(\mu)F_{2p}(\mu)-F_{2p}(\mu)F_{2p-1}(\mu)}{F_{2p}(\mu)^2+F_{2p-1}(\mu)^2}}\\
&\geq \log\sqrt{\ttmatrix{F_{2p}(\mu)^2}{0}{0}{F_{2p}(\mu)^2}}
= \log (F_{2p}(\mu))\\
&\geq (2p-1)\log\mu+\log\l(\frac{9}{10}\r),
\end{aligned}
$$
where we used Lemma \ref{lm:Fbounds} (ii) and \eqref{eq:x0est} again in the last step. The same estimates apply to $\log|(B^{2p})^\dagger|$. Upon returning to \eqref{eq:conclusion}, we find
$$
\begin{aligned}
\gam_1
\geq&
\goth{p} \l((2p-1)\log\mu+\log\l(\frac{9}{10}\r)\r)\\
&-(1-\goth{p})\l(p\log \mu+\log\l(\frac{160}{9\de_1}\r)\r)
-10^{11}\mu^{-2p}\de_2^{-2}\max\{1,(\de_1 \de_2)^{-2}\}\\
\geq& \frac{1}{2}\goth{p}p\log\mu,
\end{aligned}
$$
where we liberally used $\goth{p}\geq \frac{1}{2}$ and $\mu\geq b-1\geq b_0-1$ together with \eqref{eq:bass} in the last step. This proves Theorem \ref{thm:schro}.
\e{proof}


\section{Quantitative stability results for the Lyapunov exponent}
\label{sect:pert}
In this section, we utilize the effective Avalanche Principle to obtain quantitative stability results for Lyapunov exponents near sequences of aligned diagonal matrices (Proposition \ref{prop:stability1} and Theorem \ref{thm:stability2}). In general, continuity results for the Lyapunov exponent (and their limitations) are a topic of significant ongoing interest in the dynamical systems community; see for example \cite{BoVi,bourgain_book,DK1,DK2,DKS,GS,Ruelle,Schlag,viana_book,viana18} and references therein. The results presented here differ, as far as we can see, from the existing ones in two main ways:
\be{enumerate}

\item
We prove quantitative, deterministic continuity results using the norm topology of the matrices. The bounds are completely independent of a potentially underlying cocycle structure (both algebraically and dynamically).

\item The results require perturbing around at least slightly aligned matrices. With regards to applications to Schr\"odinger cocycles, this means other considerations, i.e., by Combes-Thomas estimates and Johnson's theorem \cite{Johnson,Zhang} also imply positivity of the Lyapunov exponent. Of course, the availability of these spectral-theoretic bounds is not perturbatively stable, in contrast to our bound. Moreover, even when restricting to the Schr\"odinger framework, these methods give the precise asymptotic growth rate of the Lyapunov exponent off the spectrum \eqref{eq:asymptotic} rather directly and without any additional input (e.g., without ever using Combes-Thomas estimates).  \e{enumerate}

The key technical ingredient of this section is the effective Avalanche Principle, and the idea of gap amplification by blocking which appears in the proof of Theorem \ref{thm:stability2}. This section does not use the $n$-matrix Golden-Thompson inequality.


Our investigation here is partly motivated by questions raised in the works of Chapman-Stolz \cite{CS} and  Duarte-Klein \cite{DK3} in the context of random Jacobi operators and we briefly consider that scenario in Corollary \ref{cor:silvius}.



\subsection{Stability of the Lyapunov exponent near sequences with a strongly expanding direction}

We describe the large-$x$ asymptotics of the Lyapunov exponent for matrices of the form
\beq\label{eq:xPkdefn}
L_k=r_kP+M_k \, ,
\eeq
where $P$ is a fixed rank-$1$ projection, the $\{M_k\}_{k\in\N}$ are uniformly bounded with bounded inverse, and $\{r_k\}_{k\in\N}\subset \R$ are sufficiently large parameters. 
Let $0<C_0<C_1$. We define $\curly{M}(C_0,C_1)$ to be the set of matrices $M$ such that $\|M\|\leq C_0$ and $\|M^{-1}\|\leq C_1$.

\be{prop}[Quantitative Stability Result 1]\label{prop:stability1}
Let $0<C_0<C_1$. There exists $r_0(C_0,C_1)>0$ and $C_2(C_0,C_1)>0$ so that for all $\{r_k\}_{k\in\N}$ with $|r_k|\geq r_0$ the following holds. For every $\{M_k\}_{k\in [n]}\subset \curly{M}(C_0,C_1)$ and every rank-$1$ projection $P$ define $L_k(x)$ by \eqref{eq:xPkdefn}. Then

\beq\label{eq:stability1}
\l |\frac{1}{n}\log \l\|\prod_{k=n}^1 L_k\r\|
-
\frac{1}{n}\sum_{k=1}^n \log|r_k|
\r |
\leq C_2\l( r_0^{-1/4}+2\frac{\log|r_1 r_n|}{n}\r)
\eeq
\e{prop}

Proposition \ref{prop:stability1} is proved in Section \ref{ssect:pfpropstability}.

\be{rmk}\
\be{enumerate}[label=(\roman*)] 
\item In practice, the $\{r_k\}_{k\in \N}$ are large, but $n$-independent absolute value and so upon sending $n\to\infty$ in \eqref{eq:stability1} the last term disappears.
\item Note that $\frac{1}{n}\sum_{k=1}^n \log|r_k|=\frac{1}{n}\log \|\prod_{k=n}^1 r_k P\|$ is the Lyapunov exponent of the main term, so Proposition \ref{prop:stability1} is indeed a stability result.
\item The dependence of $r_0,C_2$ on the parameters $C_0,C_1$ can be made explicit and the decay rate $1/4$ of the error term $r_0^{-1/4}$ can be replaced by any number $<1/2$.   
\e{enumerate}
\e{rmk}

\subsection{Application to Schr\"odinger cocycles: Asymptotics off the spectrum}
\label{ssect:asympoffspec}

In this section, we use the framework of Proposition \ref{prop:stability1} to study
the asymptotics of the Lyapunov exponent for Schr\"odinger cocycles off the spectrum. Let $(X,\mu,T)$ be an ergodic dynamical system and let $v:X\to\R$ be continuous. As in \eqref{eq:schrodefn}, we consider the cocycle defined by the map $A_E:X\to \mathrm{GL}_d(\R)$,
$$
A_E(x):=\ttmatrix{E-v(x)}{-1}{1}{0},\qquad x\in X
$$
and we denote $L_k(x):=A_E(T^k x)=r_kP+M$ with
$$
r_k=E- v(T^k x),\qquad P=\ttmatrix{1}{0}{0}{0}, \qquad \textnormal{ and }\qquad  M=\ttmatrix{0}{-1}{1}{0}\, .
$$
In this setting, Proposition \ref{prop:stability1} can be refined and the relevant constants can be computed exactly. We recall that $\mathrm{ran}(v)\subset \R$ denotes the range of the function $v$.

\be{cor}[Asymptotics of the Lyapunov exponent off the spectrum]
\label{cor:CT}
Let $r_0\geq 32$ and assume that $\mathrm{dist}(E,\mathrm{ran}(v))\geq r_0$. Then
\beq\label{eq:corCT}
\l |\frac{1}{n}\log \l\|\prod_{k=n}^1 L_k(x)\r\|
-
\frac{1}{n}\sum_{k=1}^n \log|E-v(T^k x)|
\r |
\leq 
\frac{\log|E-v(T^n x)||E-v(x)|}{n}+
1.2*10^3 r_0^{-2} \, .
\eeq
\e{cor}

For a fixed energy $E$ and $\mathrm{ran}(v)$ compact, we also have an upper bound on $\sup_{y\in X}|E-v(y)|$ and so the first term on the right-hand side of \eqref{eq:corCT} is an error term that vanishes in the $n\to\infty$ limit.

By the ergodic theorem and continuity, the correct asymptotic for the Lyapunov exponent is thus given by
\beq\label{eq:asymptotic}
\lim\limits_{n\to\infty}\frac{1}{n}\sum_{k=1}^n \log|E-v(T^k x)|=\int \log |E-v(y)|\d\mu(y),\qquad \textnormal{for } \mu-a.e.\ x \, .
\eeq

As described in Section \ref{sect:schrodinger}, the matrices $L_k$ arise as the transfer matrices of discrete Schr\"odinger operators \eqref{eq:Hdefn} whose spectrum satisfies the containment
$
\mathrm{spec}\, H(x)\subset[-2,2]+\mathrm{ran}(v).
$
Therefore, the assumption $\mathrm{dist}(E,\mathrm{ran} (v))\geq r_0\geq 32$ implies that $E\not\in \mathrm{spec}\, H(x)$.

The logarithmic growth of the Lyapunov exponent in the distance to the spectrum that follows from these considerations and Corollary \ref{cor:CT} is a well-known consequence of the Combes-Thomas estimate \cite{CT,K}. Corollary \ref{cor:CT} has the advantages of giving a precise asymptotic for the Lyapunov exponent and of being stable under perturbations of the matrices, e.g., perturbations leading to a loss of the $\mathrm{SL}(2)$ structure and thus the connection to a Schr\"odinger operator. It has the disadvantage of only applying well away from the spectrum, since it requires the distance to the spectrum be sufficiently large, namely exceeding $r_0+2\geq 34$, while the Combes-Thomas approach \cite{CT,K} yields positive lower bounds on the Lyapunov exponent immediately off the spectrum.\footnote{We mention that there is some degree of flexibility in the choice of $r_0$, since one can use alternative versions of the effective Avalanche Principle, e.g., Theorem \ref{thm:AP2}.}

\subsection{Continuity of the Lyapunov exponent near aligned diagonal matrices}
For Proposition \ref{prop:stability1}, the validity of the effective Avalanche Principle is ensured by the large parameter $x$. Using the ``gap amplification by blocking'' insight from Section \ref{sect:schrodinger}, we can also prove a continuity result which applies to diagonal matrices which do not feature a large parameter.

\be{defn}[Aligned diagonal matrices]
Given numbers $\Gam\in (0,1)$ and $\eta,C_0,C_1>0$ we define $\curly{D}(\Gam,\eta,C_0,C_1)$ as the set of diagonal matrices $D=\mathrm{diag}(\lam_1^{(k)},\lam_2^{(k)},\ldots,\lam_d^{(k)})$ with $\lam_j^{(k)}\in\R$ that are 
\be{itemize}
\item[(i)] uniformly gapped, i.e., $\frac{|\lam^{(k)}_1|}{ \max_{2\leq i\leq d} |\lam_i^{(k)}|} \geq \Gam^{-1}$ for all $k\in \N$.
\item[(ii)] uniformly bounded and uniformly bounded away from zero, i.e.,  $\min_{2\leq i\leq d} |\lam_i^{(k)}|\geq \eta$ and $C_0\leq |\lam^{(k)}_1|\leq C_1$ for all $k\in \N$.
\e{itemize}
\e{defn}

Note that $C_0\geq \eta$, so one can replace the constant $C_0$ by $\eta$ everywhere below if one is not interested in the detailed parameter dependencies.

\be{thm}[Quantitative Stability Result 2]\label{thm:stability2}
Let $\Gam\in (0,1)$ and $\eta,C_0,C_1>0$.
For every $\eps_1>0$, there exists $\de_0(\eps_1,\Gam,\eta,C_0,C_1)\in (0,1]$ such that the following holds for all $\de\in(0,\de_0)$. For every sequence of diagonal matrices $\{D_k\}_{k\in \N}\subset \curly{D}(\Gam,\eta,C_0,C_1)$ and every sequence of matrices $\{M_k\}_{k\in N}$ satisfying
$$
\|D_k-M_k\|\leq \de,\qquad \forall k\in \N,
$$
we have
\beq\label{eq:thmstability2}
\l |\frac{1}{n}\log \l\|\prod_{k=n}^1 M_k\r\|
-
\frac{1}{n}\log \l\|\prod_{k=n}^1 D_k\r\|
\r |
\leq 
\eps_1+\frac{\nu}{n}\log (3C_1)\qquad \forall n\in \N\,,
\eeq
with
\beq\label{eq:nudefn}
\nu:=\ceil*{\frac{\log (4000/\eps_1)}{\log(\Gam^{-1})}}\in \N\,.
\eeq
\e{thm} 
Theorem \ref{thm:stability2} is proved in Section \ref{ssect:pfthmstability}. 
\be{rmk}\ 
\be{enumerate}[label=(\roman*)] 
\item The constants $\eps_1,\de_0$ are independent of $n$. To obtain a continuity result for Lyapunov exponents (i.e., infinite products), one fixes $\de>0$ and sends $n\to\infty$ so that the final term in \eqref{eq:thmstability2} disappears. 
\item The entire statement is deterministic and can be applied to ergodic cocycles if the assumptions hold almost-surely with respect to the ergodic measure $\mu$.
\item For the sequence of diagonal matrices $D_k$, one can compute the ``Lyapunov exponent'' exactly, i.e.,
$$
\frac{1}{n}\log 
\l\|\prod_{k=n}^1 D_k\r\|
=\frac{1}{n}\log \l |\prod_{k=n}^1 \lam_1^{(k)}\r |
=\frac{1}{n}\sum_{k=1}^n\log |\lam_{1}^{(k)}|
$$
for all $n\geq 1$. 


\item The argument is quantitative and the size of the allowed perturbation $\de_0(\eps_1,\Gam,\eta,C_0,C_1)$ can be made completely explicit. Specifically, defining the constants
$$
C_2:=\nu 2^\nu C_1^\nu\qquad \text{and} \qquad C_3:=\frac{2^{\nu}C_1^\nu C_2}{\eta^\nu}\, ,
$$
one can take
\beq\label{eq:delta0ass}
\de_0:=\min\l\{1,\frac{\eta}{2},\frac{\eta^\nu}{2C_2},\frac{C_0^{2\nu}}{2^{\nu+1} C_1^\nu C_2},\frac{\Gam^{-\nu}C_1^\nu C_3}{3},\frac{\eps_1}{10(C_2+C_2)^2}\r\}.
\eeq

\e{enumerate}
\e{rmk}

\subsection{Application of Theorem \ref{thm:stability2} to Lyapunov exponents of Jacobi operators}
The following models were considered in \cite{CS,DK3}.
Let $E\in \R$ and $\{\theta_j\}_{j\in \Z}\subset (0,\infty)$. We consider the matrices
$$
M_j(E)=\ttmatrix{\frac{E^2-1}{\theta_j}}{-E\theta_j}{\frac{E}{\theta_j}}{-\theta_j}.
$$
The interest in these arises because they are the two-step transfer matrices of the Jacobi operator $J:\ell^2(\Z)\to\ell^2(\Z)$ defined by
$$
(J\psi)_n
=\be{cases}
\theta_{n} \psi_{n-1}+\psi_{n+1},\qquad  &\textnormal{if $n$ is even},\\
\psi_{n-1}+\theta_n \psi_{n-1}, \qquad  &\textnormal{if $n$ is odd}.
\e{cases}
$$
In applications, we think of an ergodic environment obtained by sampling along orbits of a dynamical system. This can be formulated as a cocycle as described in the introduction by letting $T:X\to X$ be an ergodic invertible dynamical system which is sampled along a continuous function  $f:X\to \R$ by setting $\theta_j(x)=f(T^j x)$ for all $j\in\Z$ and all $x\in X$.

One question of interest in \cite{CS,DK3} is the stability of the Lyapunov exponent at $E=0$. The following corollary of the deterministic Theorem \ref{thm:stability2} says that the Lyapunov exponent is stable at $E=0$ if the $|\theta_j|$ are uniformly larger or uniformly smaller than $1$. 

\be{cor}[Stability at $E=0$]\label{cor:silvius}
Suppose that either
\beq\label{eq:assf}
\begin{aligned}
\mathrm{supp}\{\theta_j\}_{j\in \Z}\i\subset (0,1)
\qquad \textnormal{ or }\qquad 
\mathrm{supp}\{\theta_j\}_{j\in \Z}\i\subset (1,\infty) \, .
\end{aligned}
\eeq
Then, for every $\eps_1>0$, there exists 
$E_0(\eps_1,\beta_1,\beta_2)>0, \xi(\eps_1,\beta_1,\beta_2)>0$ such that $
|E|<E_0
$
implies 
\beq
\l |\frac{1}{n}\log \l\|\prod_{k=n}^1 M_k(E)\r\|
-
\frac{1}{n}\log \l\|\prod_{k=n}^1 M_k(0)\r\|
\r |
\leq 
\eps_1+\frac{\xi}{n} \, .
\eeq
\e{cor}

We note that $E=0$ is not in the spectrum of $J$ when \eqref{eq:assf} holds.  Extending this stability result to more general distributions of $\{\theta_j\}$ for which $E=0$ is in the spectrum is an interesting open problem. This likely requires relaxing the deterministic assumption on the Avalanche Principle. 

\be{proof}
By condition \eqref{eq:assf}, there exist $\beta_1,\beta_2>0$ such that either $\beta_1\leq \theta_j \leq 1-\beta_2$ $\forall j \in \Z$ or $1+\beta_1\leq \theta_j \leq \beta_2$ $\forall j \in \Z$. Let $\eps_1>0$. We apply Theorem \ref{thm:stability2} to the diagonal matrices.
$$
D_k=M_k(0)=\ttmatrix{\theta_j}{0}{0}{\frac{1}{\theta_j}} \, .
$$
Condition \eqref{eq:assf} ensures that the assumptions of Theorem \ref{thm:stability2} on the $\{D_k\}$ are satisfied for appropriate $\Gam,\eta,C_1,C_2$ depending only on $\beta_1,\beta_2$. Hence, Theorem \ref{thm:stability2} yields the existence of a $\de_0(\eps_1,\beta_0,\beta_1)>0$. Finally, note that
$$
\|M_k(E)-M_k(0)\|
\leq 
\l\|
\ttmatrix{0}{\frac{E}{\theta_j}}{E\theta_j}{\frac{-E^2}{\theta_j}}
\r\|
\leq C(\beta_1,\beta_2) |E| \, ,
$$
where $C(\beta_1,\beta_2)>0$ is an appropriate constant. Corollary \ref{cor:silvius} now follows by setting $E_0:=\de_0/C(\beta_1,\beta_2)>0$.
\e{proof}

\be{rmk}More general block Jacobi matrices to which the stability results can be applied appear in \cite{CS}. Here we only focus on the special case where the $M_j$ are $2\times 2$ for simplicity as in \cite{DK3}. 
\e{rmk}

\subsection{Proofs of the Stability results from the Avalanche Principle}
In this section, we give the proofs of Proposition \ref{prop:stability1}, Corollary \ref{cor:CT}, and Theorem \ref{thm:stability2}.

\subsubsection{Proof of Proposition \ref{prop:stability1}}
\label{ssect:pfpropstability}
\be{lm}\label{lm:APverifx}
There exists $r_0(C_0,C_1)>1$ so that for all $\{r_k\}_{k\in \N}\subset (r_0,\infty)$ the sequence $\{L_k\}_{k\in \N}$ satisfies the conditions of $\mathbf{AP}\l(\frac{1}{10},\frac{1}{10},5,11\r)$ with $\eps=\frac{1}{10}$ and $\kappa=\frac{1}{1000 r_0^{1/4}}$.
\e{lm}

\be{proof}[Proof of Lemma \ref{lm:APverifx}]
We have
$$
L_{k}^\dagger L_k=r_k^2 P+r_kM_{k}P+r_kPM_k+M_k^2 \, .
$$
Hence, by Weyl's Theorem on norm perturbations of Hermitian matrices and $M_k\in\curly{M}(C_0,C_1)$, there exists a constant $C>0$ such that
\beq\label{eq:lmAP1}
|s_1(L_k)^2-r_k^2|\leq C|r_k| \qquad \text{and} \qquad
s_2(L_k)^2\leq C|r_k| \, .
\eeq
The former implies that $|s_1(L_k)-|x||\leq C$ and so
$$
\mathrm{gr}(L_k)=\frac{s_1(L_k)}{s_2(L_k)}\geq \frac{|r_k|-C}{C'|r_k|^{1/2}}\geq C r_0^{1/2}\geq 1000\, r_0^{1/4}
$$
for sufficiently large $r_0$. This verifies the alignment assumption (G).

A similar argument based again on Weyl's Theorem for norm perturbations and $M_k\in\curly{M}(C_0,C_1)$ yields
\beq\label{eq:lmAP2}
|\|L_{k+1}L_k\|-|r_{k+1}r_k||\leq C  \max\{|r_k|,|r_{k+1}|\}
\eeq
for large enough $x$.  Therefore
$$
\rho(L_{k},L_{k+1})=\frac{\|L_{k+1}L_k\|}{\|L_{k+1}\|\|L_k\|}
\geq\frac{|r_{k+1}r_k|-C\max\{|r_k|,|r_{k+1}|\}}{(|r_k|+C')(|r_{k+1}|+C'')}\geq \frac{1}{10}
$$
for sufficiently large $r_0$. This verifies the alignment assumption (A) and proves Lemma \ref{lm:APverifx}.
\e{proof}
\be{proof}[Proof of Proposition \ref{prop:stability1}]
By Lemma \ref{lm:APverifx}, we can apply $\mathbf{AP}\l(\frac{1}{10},\frac{1}{10},5,11\r)$ with $\eps=\frac{1}{10}$ and $\kappa=\frac{1}{1000 r_0^{1/4}}$. Combining this with  \eqref{eq:lmAP1} and \eqref{eq:lmAP2} we obtain the upper bound
$$
\begin{aligned}
\frac{1}{n}\log \l\|\prod_{k=n}^1 L_k\r\|
&\leq \frac{1}{n}\sum_{k=1}^{n-1}\log \l\|L_{k+1}L_k\r\|-\frac{1}{n}\sum_{k=1}^{n}\log \l\|L_k\r\|+\frac{11\kappa}{\eps^2}\\
&\leq \frac{1}{n}\sum_{k=1}^{n-1} \log
\l(\frac{|r_{k+1}r_k|+C\max\{|r_k|,|r_{k+1}|\}}
{|r_k|-C|r_k|^{1/2}}
\r)
-\frac{1}{2n}\log|r_n|+Cr_0^{-1/4}\\
&\leq \frac{1}{n}\sum_{k=1}^{n}\log|r_k| +\frac{1}{n}\sum_{k=1}^{n-1} \log
\l(\frac{1+C\max\{|r_k|^{-1},|r_{k+1}|^{-1}\}}
{1-C|r_k|^{-1/2}}
\r)
+Cr_0^{-1/4}\\
&\leq \frac{1}{n}\sum_{k=1}^{n} \log |r_k|+Cr_0^{-1/4}
\end{aligned}
$$
for sufficiently large $r_0$.  An analogous argument based on the lower bound in $\mathbf{AP}\l(\frac{1}{10},\frac{1}{10},5,11\r)$ yields
$$
\frac{1}{n}\log \l\|\prod_{k=n}^1 L_k\r\|\geq \frac{1}{n}\sum_{k=1}^{n} \log |r_k|-Cr_0^{-1/4}-2\frac{\log|r_1|}{n}-2\frac{\log|r_n|}{n}
$$
for sufficiently large $r_0$. This proves Proposition \ref{prop:stability1}.
\e{proof}

\subsubsection{Proof of Corollary \ref{cor:CT}}
We denote $r_k=E-v(T^k x)$ and $L_k=L_k(x)$. The proof is similar to the one of Proposition \ref{prop:stability1}, but we exhibit explicit constants.

\be{lm}\label{lm:normestimates}
The matrices $\{L_k(x)\}_{k\in\N}$ satisfy the bounds
\begin{align}
\label{eq:Lknorm}
|r_k|\leq&\|L_k\|\leq |r_k|(1+r_0^{-2})		\\
\label{eq:Lk+1knorm}
|r_{k+1}||r_k|\leq &\|L_{k+1}L_k\|
\leq |r_{k+1}||r_k| (1+3r_0^{-2})  \, .
\end{align}
\e{lm}

\be{proof}[Proof of Lemma \ref{lm:normestimates}]
 First, we have 
\begin{align}\label{eq:psidefn}
\|L_k\|^2= 1 + \frac{1}{2} \l( r_k^2 +|r_k|  \sqrt{4+r_k^2}\r)=:\psi(r_k).
\end{align}
The elementary inequality $\sqrt{1+y}\leq 1+y/2$ for all $y>0$ implies that
\beq\label{eq:psiestimate}
x^2\leq \psi(x)\leq x^2+2
\eeq
and after taking square roots $|r_k|\leq \|L_k\|\leq |r_k|+|r_k|^{-1}$ from wich \eqref{eq:Lknorm} follows.

Second, we have the exact formula 
\begin{align}
\|L_{k+1}L_k\|= \frac{1}{\sqrt{2}}\sqrt{\phi(r_k,r_{k+1}) +\sqrt{\phi(r_k,r_{k+1})^2 -4} }\, ,
\end{align}
where we introduced the function
\begin{align*}
\phi(x,y):=2-2xy+y^2+x^2(1+y^2)=2+(x-y)^2+x^2y^2
\end{align*}
which satisfies the estimates 
\beq\label{eq:phiestimate}
2+x^2y^2 \leq \phi(x,y)\leq 2+(|x|+|y|)^2+x^2y^2
\leq \l(\frac{5}{\min\{|x|^2,|y|^2\}}+1\r) x^2y^2
\eeq
for all $|x|,|y|\geq 1$. By $|r_k|,|r_{k+1}|\geq r_0\geq 32$, these imply
$$
|r_{k+1}||r_k|\leq \|L_{k+1}L_k\|
\leq  \sqrt{\phi(r_k,r_{k+1})}
\leq |r_{k+1}||r_k|\sqrt{1+5r_0^{-2}} 
\leq |r_{k+1}||r_k| (1+3r_0^{-2})
$$
and thus \eqref{eq:Lk+1knorm}. This proves Lemma \ref{lm:normestimates}.
\e{proof}

\be{proof}[Proof of Corollary \ref{cor:CT}]
We note that Lemma \ref{lm:normestimates} and $|r_k|\geq r_0\geq 32$ imply that
the matrices $\{L_k(x)\}_{k\in\N}$ satisfy the conditions of $\mathbf{AP}\l(\frac{1}{10},\frac{1}{10},5,11\r)$ with $\kappa = r_0^{-2}\leq \frac{1}{1000}$, and $\eps = 1/10$. Indeed, Lemma \ref{lm:normestimates} gives 
$$
\mathrm{gr}(L_k)=\|L_k\|^2\geq |r_k|^2\geq r_0^2\geq 1000
$$
and, using also  $|r_k|,|r_{k+1}|\geq r_0\geq 32$,
$$
\rho(L_k,L_{k+1})=\frac{\|L_{k+1} L_k\|}{\|L_{k+1}\| \|L_k\|}
\geq \frac{|r_{k+1}||r_k|}{(|r_{k+1}|+r_0^{-1})(|r_{k}|+r_0^{-1})}
\geq \frac{1}{10} \, .
$$
Thus we can apply $\mathbf{AP}\l(\frac{1}{10},\frac{1}{10},5,11\r)$ with $\kappa = r_0^{-2}\leq \frac{1}{1000}$, and $\eps = 1/10$. Combining the upper bound with Lemma \ref{lm:normestimates} gives
$$
\begin{aligned}
\frac{1}{n}\log \l\|\prod_{k=n}^1 L_k\r\|
&\leq \frac{1}{n}\sum_{k=1}^{n-1}\log \l\|L_{k+1}L_k\r\|-\frac{1}{n}\sum_{k=1}^{n}\log \l\|L_k\r\|+\frac{11\kappa}{\eps^2}\\
&\leq \frac{1}{n}\sum_{k=1}^{n-1}\log \l(|r_{k+1}||r_k| (1+3r_0^{-2})\r)
-\frac{1}{n}\sum_{k=1}^{n}\log |r_k|+1.1*10^3 r_0^{-2}\\
&\leq
\frac{1}{n}\sum_{k=1}^{n}\log |r_k|
-\frac{1}{n}\log(|r_1||r_n|)+1.2*10^3 r_0^{-2}\\
&\leq
\frac{1}{n}\sum_{k=1}^{n}\log |r_k|
+1.2*10^3 r_0^{-2} \, .
\end{aligned}
$$
Similarly, the lower bound in $\mathbf{AP}\l(\frac{1}{10},\frac{1}{10},5,11\r)$ and Lemma \ref{lm:normestimates} combine to give
$$
\begin{aligned}
&\frac{1}{n}\log \l\|\prod_{k=n}^1 L_k\r\|
\geq \frac{1}{n}\sum_{k=1}^{n-1}\log \l\|L_{k+1}L_k\r\|-\frac{1}{n}\sum_{k=1}^{n}\log \l\|L_k\r\|-\frac{5\kappa}{\eps^2}\\
&\geq \frac{1}{n}\sum_{k=1}^{n-1}\log \l(|r_{k+1}||r_k| \r)
-\frac{1}{n}\sum_{k=1}^{n}\log( |r_k|(1+r_0^{-2}))-500 r_0^{-2}\\
&\geq
\frac{1}{n}\sum_{k=1}^{n}\log |r_k|
-\frac{1}{n}\log(|r_1||r_n|)-600 r_0^{-2} \, .
\end{aligned}
$$
This finishes the proof of Corollary \ref{cor:CT}.
\e{proof}

\subsubsection{Proof of Theorem \ref{thm:stability2}}
\label{ssect:pfthmstability}
We write $D_k=\mathrm{diag}(\lam_1^{(k)},\lam_2^{(k)},\ldots,\lam_d^{(k)})$. We first note that we can assume without loss of generality that $\lam_1^{(k)}>0$ for all $k\in [n]$, by replacing $M_k$ with $-M_k$ whenever this is not the case.

Let $\eps_1\in (0,1)$. We now amplify the singular gap by blocking. For every $k\in\N$, we define a new matrix $\tilde M_k$ consisting of $\nu$-blocks of the original sequence, i.e.,
$$
\tilde M_k:=\prod_{j=k+\nu-1}^k M_j,\qquad \forall k\in \N \, ,
$$
and similarly for the diagonal matrices
$$
\tilde D_k:=\prod_{j=(k-1)\nu+1}^{k\nu} D_j,\qquad \forall k\in \N\, .
$$
These new matrix sequences satisfy the effective Avalanche Principle. 

\be{lm}\label{lm:APverif}
The sequences $\{\tilde D_k\}_{k\in \N}$ and $\{\tilde M_k\}_{k\in \N}$ satisfies the conditions of $\mathbf{AP}\l(\frac{1}{10},\frac{1}{10},5,11\r)$ with $\eps=\frac{1}{10}$ and $\kappa=\frac{\eps_1}{2000}$.
\e{lm}
We note that $\kappa/\eps^2=\frac{\eps_1}{10}\leq \frac{1}{10}$ as required.

\be{proof}[Proof of Lemma \ref{lm:APverif}]\ 

\dashuline{Step 1: Verification for $\{\tilde D_k\}_{k\in \N}$.}
The validity of the gap assumption (G) follows from the fact that the singular values of diagonal matrices are equal to the absolute value of their eigenvalues via the estimate
\beq\label{eq:tildeLkgap}
\mathrm{gr}(\tilde D_k)
= \frac{\prod\limits_{j=(k-1)\nu+1}^{k\nu}  \lam_1^{(j)}}{\max\limits_{2\leq i\leq d} \prod\limits_{j'=(k-1)\nu+1}^{k\nu} |\lam_i^{(j')}|}
\geq \prod_{j=(k-1)\nu+1}^{k\nu} \frac{ \lam_1^{(j)}}{\max\limits_{2\leq i\leq d}|\lam_i^{(j)}|}\geq \Gam^{-\nu}>\frac{4000}{\eps_1}\geq \frac{1}{\kappa}\, ,
\eeq
where we used Definition \eqref{eq:nudefn} of $\nu$. The alignment assumption (A) holds trivially for the $\{\tilde D_k\}_{k\in\N}$ because these diagonal matrices are in fact fully aligned, i.e.,
$$
\rho(\tilde D_k,\tilde D_{k+1})=\frac{\|\tilde D_{k+1}\tilde D_k\|}{\|\tilde D_{k+1}\|\|\tilde D_k\|}=1\geq \frac{1}{10}=\eps\, .
$$
This concludes step 1.

\dashuline{Step 2:  Verification for $\{\tilde M_k\}_{k\in \N}$.} We first note that our assumptions on $D_k$ and Definition \eqref{eq:delta0ass} imply
\beq\label{eq:sdlb}
s_d(M_k)=\min_{x\in \R^d\setminus \{0\}}\frac{\|M_k x\|}{\|x\|}
\geq \min_{x\in \R^d\setminus \{0\}}\frac{\|D_k x\|}{\|x\|}-\|M_k-D_k\|
\geq \eta-\de_0
\geq \frac{\eta}{2}>0\, .
\eeq
This shows that each $M_k$ is invertible and therefore also each $\tilde M_k$ is invertible. By telescoping, the triangle inequality, $\|AB\|\leq \|A\|\|B\|$, and our assumptions,
$$
\begin{aligned}
\|\tilde D_k-\tilde M_k\|
&=\l\|\prod_{j=(k-1)\nu+1}^{k\nu}  D_j-\prod_{j=(k-1)\nu+1}^{k\nu}  M_j\r\|\\
&\leq \sum_{J=(k-1)\nu+1}^{k\nu}  \l\|
\prod_{j=(k-1)\nu+1}^J D_j \prod_{j=J-1}^{k\nu} M_j
-\prod_{j=(k-1)\nu+1}^{i-1} D_j \prod_{j=J}^{k\nu} M_j\r\|\\
&\leq   \sum_{J=(k-1)\nu+1}^{k\nu} 
\l(\prod_{j=(k-1)\nu+1}^{J-1} \|D_j\|\r) \|D_J-M_J\| \l(\prod_{j=J-1}^{k\nu} \|M_j\|\r)\\
&\leq  \de   \sum_{J=(k-1)\nu+1}^{k\nu} 
\l(\prod_{j=(k-1)\nu+1}^{J-1} C_1\r) \l(\prod_{j=J-1}^{k\nu} (1+\de_0)C_1\r)\\
&\leq \de \nu 2^\nu C_1^\nu \, ,
\end{aligned}
$$
where the last estimate uses that $\de_0\leq 1$. In summary,
\beq\label{eq:tildepert}
\|\tilde D_k-\tilde M_k\|\leq \de C_2,\qquad \textnormal{for} \qquad C_2=\nu 2^\nu C_1^\nu.
\eeq
It remains to notice that the conditions (G) and (A) involve spectral data and $n$-independent constants and are therefore stable under perturbations of sufficiently small norm by Weyl's theorem. We proceed with the details, also with the goal in mind of explicitly verifying the effective choice \eqref{eq:delta0ass} for $\de_0$.

We use \eqref{eq:tildepert} to verify that the sequence $\{\tilde M_k\}_{k\in\N}$ satisfies assumptions (G) and (A) of $\mathbf{AP}\l(\frac{1}{10},\frac{1}{10},5,11\r)$ with $\eps=\frac{1}{10}$ and $\kappa=\frac{\eps_1}{2000}$. We first verify assumption (A). By the triangle inequality, our assumptions on $\{D_k\}$, and $\de_0\leq 1$, we have
\beq\label{eq:Mprodlb}
\begin{aligned}
\l|
\|\tilde M_{k+1}\tilde M_{k}\|
-\|\tilde D_{k+1}\tilde D_{k}\|
\r|
&\leq 
\|\tilde D_{k+1}\| \|\tilde D_k-\tilde M_k\|
+\|\tilde D_{k+1}-\tilde M_{k+1}\| \|\tilde M_k\|\\
&\leq \de C_1^\nu C_2
+\de (1+\de_0)^\nu C_1^\nu C_2\\
&\leq \de 2^{\nu+1} C_1^\nu C_2 \, .
\end{aligned}
\eeq
Hence, employing various estimates following from Definition \eqref{eq:delta0ass} of $\de_0$, we find that
$$
\begin{aligned}
\rho(\tilde M_k,\tilde M_{k+1})
&=\frac{\|\tilde M_{k+1}\tilde M_{k}\|}{\|\tilde M_{k+1}\|\|\tilde M_{k}\|}\\
&\geq \frac{\|\tilde D_{k+1}\tilde D_{k}\|-\de 2^{\nu+1} C_1^\nu C_2}{(\|\tilde D_{k+1}\|+\de C_2)(\|\tilde D_{k}\|+\de C_2)}\\
&\geq  \frac{\prod\limits_{j=k\nu+1}^{(k+1)\nu} \lam_d^{(j)}\prod\limits_{j'=(k-1)\nu+1}^{k\nu} \lam_d^{(j')}
-\de 2^{\nu+1} C_1^\nu C_2}{\prod\limits_{j''=k\nu+1}^{(k+1)\nu} \lam_d^{(j'')}\prod\limits_{j'''=(k-1)\nu+1}^{k\nu} \lam_d^{(j''')}+2\de C_1^\nu C_2+\de^2 C_2^2}\\
&\geq  \frac{1}{2}\,\frac{1}{1+\frac{2\de C_1 C_2+\de^2 C_2^2}{C_0^{2\nu}}}\\
&\geq \frac{1}{10}=\eps\, ,
\end{aligned}
$$ 
as desired.

 Finally, we verify assumption (G) for the $\{\tilde M_k\}_{k\in\N}$. This involves the singular values of the $\tilde M_k$, i.e., the square roots of the eigenvalues of $\tilde M_k^\dagger \tilde M_k$. By \eqref{eq:tildepert}, $\|A^\dagger\|=\|A\|$, and our assumptions on $\tilde D_k$, we have
\beq\label{eq:MDsquarediff}
\|\tilde M_k^\dagger M_k-\tilde D_k^\dagger D_k\|
\leq \|\tilde M_k\|  \|\tilde M_k- \tilde D_k\|+\|\tilde M_k-\tilde D_k\|\|\tilde D_k\|
\leq \de 2^{\nu+1}C_1^\nu C_2 \, .
\eeq
By Weyl's theorem on norm perturbations of Hermitian matrices, this implies that
$$
|s_i(\tilde M_k)^2-s_i(\tilde D_k)^2|\leq \de 2^{\nu+1}C_1^\nu C_2\qquad \forall 1\leq i\leq d,
$$
and so
$$
|s_i(\tilde M_k)-s_i(\tilde D_k)|\leq \de \frac{2^{\nu+1}C_1^\nu C_2}{s_d(\tilde M_k)+s_d(\tilde D_k)} \qquad \forall 1\leq i\leq d \, .
$$
To bound the denominator from below, we use our assumption on $D_k$ as follows,
\beq\label{eq:sdtildeDk}
s_d(\tilde D_k)
= \min_{1\leq i'\leq d}\prod\limits_{j=(k-1)\nu+1}^{k\nu}|\lam_d^{(i')}|
\geq \prod\limits_{j=(k-1)\nu+1}^{k\nu} \min_{1\leq i'\leq d}\limits|\lam_d^{(i')}|
\geq \eta^\nu \, .
\eeq
Repeating the estimates in \eqref{eq:sdlb} with $M_k, D_k$ replaced by $\tilde M_k, \tilde D_k$ and using \eqref{eq:tildepert} gives
\beq\label{eq:sdtildeMk}
s_d(\tilde M_k)\geq \eta^\nu-\de C_2\geq \frac{\eta^\nu}{2}
\eeq
and thus
\beq\label{eq:sipert}
|s_i(\tilde M_k)-s_i(\tilde D_k)|\leq \de C_3 \qquad \textnormal{for} \qquad C_3=\frac{2^{\nu}C_1^\nu C_2}{\eta^\nu} \qquad \forall 1\leq i\leq d\, .
\eeq
Combining \eqref{eq:sipert} with \eqref{eq:tildeLkgap}, our assumptions on $\tilde D_k$ and Definition \eqref{eq:delta0ass} of $\de_0$, we have
$$
\begin{aligned}
\mathrm{gr}(\tilde M_k)
&=\frac{s_1(\tilde M_k)}{s_2(\tilde M_k)}
\geq \frac{s_1(\tilde D_k)-\de C_3}{s_2(\tilde D_k)+\de C_3}
\geq \frac{\prod\limits_{j=(k-1)\nu+1}^{k\nu}  \lam_1^{(j)}-\de C_3}{\prod\limits_{j'=(k-1)\nu+1}^{k\nu}\max\limits_{2\leq i\leq d}  |\lam_i^{(j')}|+\de C_3}\\
&\geq \frac{\prod\limits_{j=(k-1)\nu+1}^{k\nu}  \lam_1^{(j')}-\de C_3}{\Gam^{\nu}\prod\limits_{j'=(k-1)\nu+1}^{k\nu} \lam_1^{(j)} +\de C_3}
\geq \Gam^{-\nu}
\frac{1-\de C_1^\nu C_3}{1+\Gam^{-\nu} \de C_1^\nu C_3}
\geq \frac{1}{2}\Gam^{-\nu}
>\frac{2000}{\eps_1}=\frac{1}{\kappa} \, ,
\end{aligned}
$$
where the last estimate holds by the Definition \eqref{eq:nudefn} of $\nu$. This proves Lemma \ref{lm:APverif}.
\e{proof}

We are now ready to give the

\be{proof}[Proof of Theorem \ref{thm:stability2}]
Writing $n=\nu\floor{n/\nu}+r$ with $0\leq r<\nu$ and using $\|M_k\|\leq 2\|D_k\|\leq  2C_1$ gives
$$
\begin{aligned} 
&\l |\frac{1}{n}\log \l\|\prod_{k=n}^1 M_k\r\|-\frac{1}{n}\log \l\|\prod_{k=n}^1 D_k\r\|
\r|\\
&\hspace{30mm}\leq \l| \frac{1}{n}\log \l\|\prod_{k=\floor{n/\nu}}^1 \tilde M_k\r\|
-\frac{1}{n}\log \l\|\prod_{k=\floor{n/\nu}}^1 \tilde D_k\r\|\r|
+\frac{\nu}{n}\log (3C_1)\\
&\hspace{30mm}\leq \frac{1}{\nu}\l| \frac{1}{\floor{n/\nu}}\log \l\|\prod_{k=\floor{n/\nu}}^1 \tilde M_k\r\|
-\frac{1}{\floor{n/\nu}}\log \l\|\prod_{k=m}^1 \tilde D_k\r\|\r|
+\frac{\nu}{n}\log (3C_1)\, .\\
\end{aligned}
$$
By Lemma \ref{lm:APverif}, the sequences $\{\tilde L_k\}_{k\in [n]}$ and $\{\tilde M_k\}_{k\in [n]}$ satisfy the conditions of  $\mathbf{AP}\l(\frac{1}{10},\frac{1}{10},5,11\r)$ with $\eps=\frac{1}{10}$ and $\kappa=\frac{\eps_1}{2000}$. Hence, setting $m=\floor{n/\nu}$, and recalling \eqref{eq:tildepert}, \eqref{eq:Mprodlb}
$$
\begin{aligned} 
\frac{1}{m}\log \l\|\prod_{k=m}^1 \tilde M_k\r\|
-\frac{1}{m}\log \l\|\prod_{k=m}^1 D_k\r\|
&\leq
\frac{1}{m}\sum_{k=1}^{m-1}\log \frac{\|\tilde M_{k+1}\tilde M_k\|}{\|\tilde D_{k+1}\tilde D_k\|}
+\frac{1}{m}\sum_{k=1}^{m}\log \frac{\|\tilde D_k\|}{\|\tilde M_k\|}+\frac{16\kappa}{\eps^2}\\
&\leq 
\frac{1}{m}\sum_{k=1}^{m-1}\log (1+\de C_2)
+\frac{1}{m}\sum_{k=1}^{m}\log (1+\de C_2^2)+\frac{16\kappa}{\eps^2}\\
&\leq \de(C_2+C_2^2)+\frac{16\kappa}{\eps^2}= \de(C_2+C_2^2)+\eps_1\frac{16}{20}\\
&\leq \eps_1\, ,
\end{aligned}
$$
where the last step holds by the final condition in \eqref{eq:delta0ass}. An analogous argument gives a lower bound $-\eps_1$ on the same quantity. This proves Theorem \ref{thm:stability2}.
\e{proof}

\section{A bound for almost commuting matrices} \label{sec_almostComm}

We note that the bounds from Theorem~\ref{thm:main} and Corollary~\ref{cor:asymptotic2} are not able to reproduce the fact that for commuting matrices the Lyapunov exponent features a simple closed from expression, see~\eqref{eq_Lyap_commute}, because the expanding directions of commuting matrices need not be aligned at all (in which case the Avalanche Principle does not apply). 
In applications, it is of interest to understand in which sense~\eqref{eq_Lyap_commute} is stable when considering a sequence of almost-commuting matrices \cite{CS,DK3}. This is the main purpose of this section. Due to technical reasons we assume in this section that all the matrices are positive definite and hence denote them by $\{A_k\}_{k \in [n]}$ (instead of $\{L_k\}_{k \in [n]}$). Furthermore, we assume that this sequence originates from a stationary and ergodic random process. For $m \leq n$ and $t \in \R$ we define 
\begin{align*}
X_{m,n}(t):=\l\| \prod_{k=m}^n A_k^{1+\ci t} \r\| \, .
\end{align*}

\noindent Very roughly, the idea of our approach is as follows:
\be{enumerate} 
\item The $n$-matrix Golden-Thompson inequality in Theorem \ref{thm_nGT} can be used to bound an appropriate average over $\log X_{1,n}(t)$ from below.
\item For almost-commuting matrices $X_{1,n}(t)\approx X_{1,n}(0)$, where the error term is proportional to the largest commutator $[A_j,A_k]$.
\e{enumerate} 

Unfortunately, we have not been able to prove a completely satisfactory bound in this way: we require an assumption on the convergence speed of the finite-size Lyapunov exponent to its limit for every fixed $t \in \R$, on which we comment further below. Still, we believe that the almost-commuting case is sufficiently natural and relevant that our modest efforts of bringing matrix analysis to bear on it deserve a brief discussion. We certainly hope that future research can shed further light on this important topic, perhaps by analyzing the assumption below more comprehensively.

Let us now present the precise assumption. First, note that the Kesten-Furstenberg theorem~\cite{furstenberg1960} implies the existence of the $t$-dependent Lyapunov exponents
\begin{align*}
\gamma_1(t)= \lim_{n \to \infty} \frac{1}{n} \log X_{1,n}(t) = \lim_{n \to \infty} \frac{1}{n} \E \log X_{1,n}(t) \, ,
\end{align*}
for every $t\in \R$, where $\gamma_1(0)=\gamma_1$.

The following assumption states that the convergence speed to $\gam_1(t)$ is of order $1/n$.

\begin{ass} \label{ass_commutative}
There exists a constant $c<\infty$ such that for all  $k\leq n$ and $t \in \R$
\begin{align*} 
\gamma_1(t) - \frac{c}{n-k+1} \leq \frac{1}{n-k+1} \log  X_{k,n}(t) \leq \gamma_1(t) + \frac{c}{n-k+1}  \quad \pr\text{-a.s.}\, .
\end{align*}
\end{ass}
 
 We remark on the feasibility of this assumption below. Under Assumption~\ref{ass_commutative} we can complete the two-step strategy mentioned above and derive a lower bound on the Lyapunov exponent for almost-commuting matrices.

The key is to handle the $t$-dependence as follows.

\begin{prop} \label{prop_almostComm}
Suppose the sequence $\{A_k\}_{k \in \N}$ satisfies Assumption~\ref{ass_commutative}. Then we have for all $t\in \R$
\begin{align*}
|\gamma_1- \gamma_1(t)|  \leq  \max \left \lbrace \sqrt{4 \ee^{5c} \max_{j,k \in \N} \l\| [A_j^{\ci t},A_k] \r\|}, 4 \ee^{5c} \max_{j,k \in \N} \l\| [A_j^{\ci t},A_k] \r\|,  \right \rbrace
\end{align*}
\end{prop}

We can use Proposition~\ref{prop_almostComm} to derive a lower bound on $\gamma_1$ that is tight in the commutative case.

\begin{cor} \label{cor_almostComm}
Suppose the sequence $\{A_k\}_{k \in \N}$ satisfies Assumption~\ref{ass_commutative}. Then
\begin{align*} 
\gamma_1 \geq \lambda_{\max} ( \E \log A_1) - \max \left \lbrace  \sqrt{4\ee^{5c}  \max_{j,k \in \N} \norm{[\log A_j, A_k]}},4\ee^{5c}  \max_{j,k \in \N} \norm{[\log A_j, A_k]} \right \rbrace \, .
\end{align*}
\end{cor}
We note that this bound reproduces the commutative result~\eqref{eq_Lyap_commute} as in this case the term with the commutator above vanishes, i.e., $[\log A_j, A_k]=0$. Moreover in an almost commuting scenario, i.e., a situation where $\norm{[\log A_j, A_k]}$ is small the bound presented above is expected to behave well. These proofs of Lemma~\ref{prop_almostComm} and Corollary~\ref{cor_almostComm} are deferred to the appendix. 

We finish with two remarks about Assumption \ref{ass_commutative}.
\be{rmk}\
\be{enumerate}[label=(\roman*)]
\item 
The convergence speed of order $1/n$ may appear to be impossible to realize: A natural reference point for convergence speed is the central limit theorem, and its fluctuations live on the much larger scale $1/\sqrt{n}$. (In fact, in the case of i.i.d.\ diagonal random matrices, the Lyapunov exponent reduces to a sum of independent random variables and one is indeed in the CLT setting.) However, there exist important examples of cocycles over much more rigid dynamical systems, namely quasi-periodic ones, where the convergence speed is of order $1/n$ on the level of the averages; see specifically Proposition 2 in \cite{Schlag}. In this context, we also mention that verifying Assumption \ref{ass_commutative} will likely be aided by using a multiscale scheme where each $A_j$ is itself given as a shorter product of transfer matrices on the previous scale. In that setting, large deviation estimates (see e.g.\ \cite{bourgain_book,GS,Schlag}) can conceivably be used to establish proximity to the limit $\gam_1(t)$ outside of certain small-measure sets of initial conditions. 
\item As one can glean from Corollary \ref{cor_almostComm}, the usefulness of our approach depends on the smallness of $e^{5c}$ times a commutator norm $\|[\log A_j,A_k]\|$. We have elected to present this from the arguably most natural perspective that $c$ is a universal constant and $\|[\log A_j,A_k]\|$ is small; other secenarios can be treated by the same method.

\e{enumerate}
\e{rmk}


\paragraph{Acknowledgements.} We thank J\"urg Fr\"ohlich, Silvius Klein, and Jeffrey Schenker for helpful discussions.
DS acknowledges support from the Swiss National Science Foundation via the NCCR QSIT as well as project No.~200020\_165843.
\begin{appendix}
\section{Proof of the Avalanche Principle, version 2} \label{ap_AP}
\be{proof}[Proof of Theorem \ref{thm:AP2}]
The proof follows the general steps in \cite{HLS}. We only summarize the necessary changes in the argument, which amount to various refined estimates and a different choice of the scaling parameter $t$ (see below).

We define $\eps_0=\frac{1}{5}$ and $c_0=\frac{1}{6}$. Using that $(\eps_0+\sqrt{1-\eps_0^2})\pi/2\leq 1.86$, the estimate $\frac{\sqrt{2}\pi}{2}\kappa \eps^{-2}$ on the Lipschitz constant in Corollary 5.7(b) in \cite{HLS} can be replaced by
\beq
1.86*\kappa \eps^{-2}.
\eeq
This replacement then propagates through the argument in \cite{HLS}, specifically to Lemma 5.14, Corollary 5.15 and the proof of Theorem 5.5.

Another change is made starting with Lemma 5.12 in \cite{HLS}. We define $\tilde \eps=t\eps$ with $t=\frac{72}{100}$. (This choice of $t$ is close to the largest one allowed here which is important for having a good bound in \eqref{eq:L1} below.) Lemma 5.12 in \cite{HLS} then holds true with this choice of $t$ under our present assumptions on $\kappa$ and $\eps$. To see this, we note that condition (5.7) in \cite{HLS} can be rearranged as 
\beq
t>c_0\sqrt{\frac{1+2c_0^2\eps^2}{1+c_0^2\eps^2}} \, ,
\eeq
which is comfortably satisfied for all $\eps\in (0,\eps_0)$. To establish Lemma 5.12, it then remains to check that 
\beq
f\l(\frac{c_0\eps}{t}\sqrt{1-t^2\eps^2},\sqrt{1-\frac{\eps^2}{1+2c_0^2\eps^2}}\r)<\sqrt{1-t^2\eps^2}
\eeq
for $f(x,y)=x\sqrt{1-y^2}+y\sqrt{1-x^2}$, which is tedious, but elementary. (Here we avoid estimating the square roots by their linearization as was done in (5.8) of \cite{HLS}.)  

The previously established estimates yield the following modification of Lemma 5.16 in \cite{HLS}. The relevant Lipschitz constant is now
\beq\label{eq:L1}
L=1.86*\kappa \tilde\eps^{-2}\leq 1.86*c_0 t^{-2}<\frac{3}{5}<1 \, ,
\eeq
where the fact that this is strictly less than $1$ is important for having a contraction. The relevant estimate after summing the geometric series in the proof of Lemma 5.16 then yields
\beq
\frac{\kappa \tilde\eps^{-1}}{1-L} = \frac{1/t}{1-L}\frac{\kappa}{\eps } <\frac{7}{2}\frac{\kappa}{\eps} 
\eeq
and then the right-hand side replaces the upper bound from Lemma 5.16 (which was $3\kappa/\eps$ in\cite{HLS}). 

For the proof of Theorem 5.5(i), it suffices to note that
\beq
\frac{7}{2}\frac{\kappa}{\eps} +\frac{1}{2}\leq \frac{7}{2}c_0 \eps_0+\frac{1}{2}=\frac{7}{30}+\frac{1}{2}<\sqrt{1-t^2\eps_0^2}<\sqrt{1-\tilde\eps^2}\, .
\eeq
In the conclusion of Theorem 5.5(i), the bound $3\frac{\kappa}{\eps}$ is replaced by $\frac{7}{2}\frac{\kappa}{\eps}$ due to our alternative version of Lemma 5.16 described above.

Finally, in the proof of Theorem 5.5(ii), we replace every instance of $3\frac{\kappa}{\eps}$ by $\frac{7}{2}\frac{\kappa}{\eps}$ in (5.17)-(5.19). This has the effect of replacing the lower bound in (5.19) by
\beq\label{eq:eps/3}
\frac{\eps}{\sqrt{1+2\frac{\kappa^2}{\eps^2}}}-\frac{7}{2}\frac{\kappa}{\eps}
\geq \frac{\eps}{\sqrt{1+2\frac{\kappa^2}{\eps^2}}}-\frac{7}{12}\eps
> \frac{1}{3}\eps
\eeq
and therefore the upper bound in (5.20) by $\frac{21}{2}\frac{\kappa}{\eps^2}$. Combining this with the bound (5.21) in \cite{HLS} and using $\kappa\leq c_0\eps_0^2=\frac{1}{150}$ yields the lower bound on the main quantity
\beq
\frac{\rho(L_1,L_2,\ldots,L_{n})}{\rho(L_1,L_2)\cdots \rho(L_{n-1}, L_{n})}\geq 
\ee^{-11 n\kappa/{\eps^2}}
\eeq
so $c_l=11$ as desired.

It remains to prove the upper bound on the main quantity. We first note that the upper bound $\kappa'$ in Lemma 5.17 needs to be replaced by $\kappa_j'=\l(\frac{3}{5}\r)^{j-1} \frac{6}{5}\kappa$. Indeed, Lemma 5.16 (with the constant $\frac{7}{2}\frac{\kappa}{\eps}\leq \frac{7}{12}\eps=:r$) and Proposition 2.29 in \cite{DK2} imply $\|(D\hat{g_0})_{\hat{\goth{v}}(g^j)}\|\leq \kappa\frac{r+\sqrt{1-r^2}}{1-r^2}\leq \frac{60}{53}\kappa$. From there, applying the chain rule as in \cite{HLS} implies 
\beq
\mathrm{gr}(g^j)\leq L^{j-1}\frac{60}{53}\kappa<\kappa',
\eeq
with the Lipschitz constant $L=1.86*\kappa \tilde\eps^{-2}<\frac{3}{5}$ defined in \eqref{eq:L1}. Equipped with this new version of Lemma 5.17 and recalling \eqref{eq:eps/3}, the relevant bound becomes
\beq
\log\sqrt{1+2\frac{(\kappa_j')^2}{\eps/3}}
\leq \frac{(\kappa_j')^2}{\eps/3}
\leq 7.2\l(\frac{3}{5}\r)^{j}\frac{\kappa^2}{\eps}.
\eeq
Since the new upper bound in (5.20) is $\frac{21}{2}\frac{\kappa}{\eps^2}$ and we assume $n\geq 36$, we have 
\beq
\frac{21}{2}n+7.2\sum_{j=1}^n \l(\frac{3}{5}\r)^{j}
\leq \frac{21}{2}n+7.2\frac{5}{2}
\leq 11n
\eeq
and hence the upper bound on the main quantity
\beq
\frac{\rho(L_1,L_2,\ldots,L_{n})}{\rho(L_1,L_2)\cdots \rho(L_{n-1}, L_{n})}
\leq 
\ee^{11 n\kappa/{\eps^2}}
\eeq
as desired. This proves Theorem \ref{thm:AP2}.
\e{proof}

\section{Proofs for the almost-commuting case}
\subsection{Proof of Proposition~\ref{prop_almostComm}}
Note that for each index $1\leq k\leq n$ it holds that $A_k^{1+\ci t}=A_k A_k^{\ci t}=A_k^{\ci t} A_k$. We start by rewriting the difference as a telescopic sum
$$
A_1^{1+\ci t}A_2^{1+\ci t}\ldots A_n^{1+\ci t}-A_1\ldots A_n A_1^{\ci t}\ldots A_n^{\ci t}
=\sum_{j=1}^{n-1} A_1^{1+\ci t}\ldots A_{j-1}^{1+\ci t} A_j [A_j^{\ci t},A_{j+1}\ldots A_n]A_{j+1}^{\ci t}\ldots A_n^{\ci t} \, .
$$
We express the long commutator $[A_j^{\ci t},A_{j+1}\ldots A_n]$ as a sum of individual commutators $[A_j^{\ci t},A_k ]$ by iteratively applying the Leibniz rule for commutators,
$$
[B,CD]=C[B,D]+[B,C]D \, .
$$
 We find
\begin{align*}
 [A_j^{\ci t},A_{j+1}\ldots A_n]
=\sum_{k=j+1}^n  A_{j+1}\ldots A_{k-1}  \l[A_j^{\ci t}, A_k \r] A_{k+1}\ldots A_{n}  \, .
\end{align*}
The unitary invariance of the operator norm then gives for $\eps_t:=\max_{j,k \in \N} \| [A_j^{\ci t},A_k]\|$
\begin{align*}
X_{1,n}(t) \leq X_{1,n}(0) + \eps_t \sum_{j<k\leq n} X_{1,j-1}(t) X_{j,j}(0) X_{j+1,k-1}(0) X_{k+1,n}(0) \, .
\end{align*}
Thus we obtain
\begin{align*}
 \frac{X_{1,n}(t)}{X_{1,n}(0)}
&\leq 1 + \eps_t \sum_{j<k \leq n} \frac{X_{1,j-1}(t) X_{j,j}(0) X_{j+1,k-1}(0) X_{k+1,n}(0)}{X_{1,n}(0)} \\
&\leq 2\eps_t \sum_{j<k \leq n} \frac{X_{1,j-1}(t) X_{j,j}(0) X_{j+1,k-1}(0) X_{k+1,n}(0)}{X_{1,n}(0)} \, ,
\end{align*}
where the final step is valid because we can assume without loss of generality that $ 1 \leq \eps_t \sum_{j<k \leq n} \frac{X_{1,j-1}(t) X_{j,j}(0) X_{j+1,k-1}(0) X_{k+1,n}(0)}{X_{1,n}(0)} $, as otherwise the assertion of Proposition~\ref{prop_almostComm} holds trivially.

Assumption~\ref{ass_commutative} then gives
\begin{align*}
\frac{X_{1,j-1}(t) X_{j,j}(0) X_{j+1,k-1}(0) X_{k+1,n}(0)}{X_{1,n}(0)} 
&\leq \frac{\ee^{(j-1) \gamma_1(t) + c} \ee^{\gamma_1 + c} \ee^{(k-j-1) \gamma_1 + c} \ee^{(n-k) \gamma_1+ c}}{\ee^{n \gamma_1 - c}} \\
&\leq \ee^{(j-1) (\gamma_1(t) - \gamma_1)} \ee^{5c} \, .
\end{align*}
We thus find
\begin{align*}
\frac{X_{1,n}(t)}{X_{1,n}(0)}
&\leq 1+2 \ee^{5c} \eps_t \sum_{j<k\leq n} \ee^{j (\gamma_1(t) - \gamma_1)} \\
&\leq 1+ 2\ee^{5c} \eps_t \sum_{j<k\leq n} \ee^{j |\gamma_1(t) - \gamma_1|} \\
&\leq 1+ 2\ee^{5c} \eps_t \sum_{j=0}^n (n-j) \ee^{j |\gamma_1(t) - \gamma_1|} \, .
\end{align*}
We next bounding the sum with an integral, i.e.,
\begin{align*}
\sum_{j=0}^n (n-j) \ee^{j |\gamma_1(t) - \gamma_1|} 
\leq \int_{0}^n  (n-x) \ee^{x |\gamma_1(t) - \gamma_1|} \d x + \frac{ \ee^{n |\gamma_1(t) - \gamma_1| }}{|\gamma_1(t) - \gamma_1|} \, ,
\end{align*}
which is correct because the function $\phi_{n,c}: [0,n] \ni x \mapsto (n-x)\ee^{x c}$ for $c>0$ is monotonically increasing in $[0,x^\star]$ and monotonically decreasing in $[x^\star,n]$ where $x^\star = n-1/c$. Furthermore we have $\phi_{n,c}(x^\star)=\frac{1}{c}\ee^{nc -1} \leq \frac{1}{c}\ee^{nc}$. Hence we find
\begin{align*}
\frac{X_{1,n}(t)}{X_{1,n}(0)}
& \leq 1 +2\ee^{5c} \eps_t \left( \frac{\ee^{n |\gamma_1(t) - \gamma_1|}}{|\gamma_1(t) - \gamma_1|^2} + \frac{ \ee^{n |\gamma_1(t) - \gamma_1| }}{|\gamma_1(t) - \gamma_1|} \right) \, .
\end{align*}
Together with  Assumption~\ref{ass_commutative} this implies 
\begin{align*}
|\gamma_1(t) - \gamma_1 | 
&\leq \frac{1}{n} \log X_{1,n}(t) - \frac{1}{n} \log X_{1,n}(0) +  \frac{2c}{n}  \\
&\leq \frac{1}{n} \log \left(1 +2\ee^{5c} \eps_t \Big( \frac{\ee^{n |\gamma_1(t) - \gamma_1|}}{|\gamma_1(t) - \gamma_1|^2} + \frac{ \ee^{n |\gamma_1(t) - \gamma_1| }}{|\gamma_1(t) - \gamma_1|} \Big) \right) + \frac{2c}{n} \, .
\end{align*}
As a result we obtain
\begin{align*}
\ee^{n |\gamma_1(t) - \gamma_1| }
\leq \left(1 +2\ee^{5c} \eps_t \Big( \frac{\ee^{n |\gamma_1(t) - \gamma_1|}}{|\gamma_1(t) - \gamma_1|^2} + \frac{ \ee^{n |\gamma_1(t) - \gamma_1| }}{|\gamma_1(t) - \gamma_1|} \Big) \right) \ee^{2c} \, .
\end{align*}
Since this is true for all $n\in \N$ we can conclude that
\begin{align*}
|\gamma_1(t) - \gamma_1| \leq \max\{\sqrt{ 4\ee^{5c} \eps_t}, 4\ee^{5c} \eps_t \} \, ,
\end{align*}
which proves the assertion. \qed

\subsection{Proof of Corollary~\ref{cor_almostComm}}
The n-matrix Golden-Thompson inequality from Theorem~\ref{thm_nGT} implies that
\begin{align*}
\int_{\R} f(t) \gamma_1(t) \d t
&= \lim_{n \to \infty} \frac{1}{n} \int_{\R} f(t) \log \l\| \prod_{k=1}^n A_k^{1+\ci t} \r\| \\
&\geq  \lim_{n \to \infty} \frac{1}{n}  \log \l\| \exp\left( \sum_{k=1}^n \log A_k \right) \r\| \\
&= \lim_{n \to \infty}  \lambda_{\max} \left( \frac{1}{n} \sum_{k=1}^n \log A_k \right) \\
&= \lambda_{\max}(\E \log A_1) \, ,
\end{align*}
where in the first step we swap the limit and the integral which is valid by the dominated convergence theorem.
Proposition~\ref{prop_almostComm} implies
\begin{align*}
\gamma_1 \geq \int_{\R} f(t) \gamma_1(t) \d t - \int_{\R} f(t) \max \left \lbrace \sqrt{ 4  \ee^{5c}\max_{j,k \in \N} \l\| [A_j^{\ci t},A_k] \r\|},4  \ee^{5c}\max_{j,k \in \N} \l\| [A_j^{\ci t},A_k] \r\| \right \rbrace \d t \, .
\end{align*}
We can use a well known bound~\cite[Lemma~3.8]{Sutter_book} to further simplify the bound\footnote{Alternatively it is also possible to use resolvent calculus to relate $\norm{ [A_j^{\ci t},A_k]}$ with $\norm{ [A_j,A_k]}$.} by using
\begin{align*}
\int_{\R} f(t) \sqrt{\max_{j,k \in \N} \l\| [A_j^{\ci t},A_k] \r\|} \d t
\leq \max_{j,k \in \N} \sqrt{\norm{[\log A_j, A_k]}} \int_{\R} f(t) \sqrt{|t|} \d t
\leq \max_{j,k \in \N} \sqrt{\norm{[\log A_j, A_k]}}  \, ,
\end{align*}
which then proves the assertion. \qed

\end{appendix}

\bibliographystyle{arxiv_no_month}
\bibliography{bibliofile}

\end{document}